\documentclass[reqno,12pt]{amsart}

\usepackage{amsmath,amssymb,mathtools}
\usepackage{mathrsfs,paralist,esint}
\usepackage[pagewise]{lineno}
\usepackage{setspace}
\usepackage{xcolor}
\usepackage{hyperref}
\hypersetup{urlcolor=red,citecolor=blue}
\usepackage[top=1in,bottom=1in,left=1in,right=1in]{geometry}
\usepackage[sort]{cite}

\allowdisplaybreaks[4]
\numberwithin{equation}{section}

\newtheorem{theorem}{Theorem}[section]
\newtheorem{lemma}[theorem]{Lemma}
\newtheorem{proposition}[theorem]{Proposition}
\theoremstyle{definition}

\newcommand{\dd}{\;\mathrm{\,d}}
\newcommand{\norm}[1]{\left\lVert #1\right\rVert}

\begin{document}

\title[Keller--Segel system with indirect signal production on closed surfaces]
{Critical mass for infinite-time blow-up in a Keller--Segel system with indirect signal production on closed Riemannian surfaces}

\author[Li]{Jing Li}
\address{J. Li: School of Mathematics, Southeast University, Nanjing 211189, P. R. China}
\email{230258239@seu.edu.cn}

\author[Li]{Yuxiang Li}
\address{Y. Li: School of Mathematics, Southeast University, Nanjing 211189, P. R. China}
\email{lieyx@seu.edu.cn}

\date{}


\begin{abstract}
We study a parabolic--elliptic--ODE chemotaxis system with indirect signal production on closed connected Riemannian surfaces.
For all sufficiently regular nonnegative initial data, the system admits
a unique global classical solution.
We show that the solution remains globally bounded if the total mass $M<8\pi$.
For each $M>8\pi$, we find initial data of mass $M$ that yield infinite-time blow-up solutions.
We then adapt the argument to the corresponding Neumann problem in a smooth
bounded planar domain and obtain infinite-time blow-up solutions for any mass $M>4\pi$.
At the critical mass $M=8\pi$, we show that the solutions are globally bounded on spheres equipped with the canonical spherical metric,
and we construct a solution blowing up in infinite time on the Riemann sphere equipped with a special metric.
\end{abstract}
\maketitle

\section{Introduction}
\label{sec:main-results}
Let $(\Sigma,g)$ be a closed connected Riemannian surface endowed with the Riemannian measure $\mu_g$.
We consider the following parabolic--elliptic--ODE chemotaxis system with indirect signal production: 
\begin{equation}\label{eq:system}
\begin{cases}
\partial_t u=\Delta_g u-\operatorname{div}_g(u\nabla_g v),
    & x\in\Sigma,\ t>0,\\[1mm]
0=\Delta_g v-\overline w+w,
    & x\in\Sigma,\ t>0,\\[1mm]
\partial_t w=-w+u,
    & x\in\Sigma,\ t>0,\\[1mm]
u(\cdot,0)=u_0,\quad w(\cdot,0)=w_0,
    & x\in\Sigma,
\end{cases}
\end{equation}
where $\overline w$ denotes the spatial average of $w$, defined by
\begin{equation}\label{eq:kdef}
\overline w(t):=\frac{1}{\mu_g(\Sigma)}\int_{\Sigma}w(x,t)\dd\mu_g,
\end{equation}
and $v$ is normalized by
\begin{equation}\label{eq:normalization}
 \int_{\Sigma}v(x,t)\dd\mu_g=0\quad \text{for all }t\ge0.
\end{equation}
The first equation in \eqref{eq:system} implies conservation of mass, namely,
\begin{equation}\label{eq:mass}
 \int_{\Sigma}u(x,t)\dd\mu_g
 =\int_{\Sigma}u_0\dd\mu_g=:M
 \quad \text{for all } t>0.
\end{equation}
Here $\nabla_g$, $\operatorname{div}_g$, and $\Delta_g$ denote the gradient, divergence, and Laplace--Beltrami operator associated with $g$, respectively.

\subsection{Direct signal production in planar domains}
The original Keller--Segel system describes how cells move towards higher concentrations of a chemical signal \cite{1970-JTB-KellerSegel}.
J\"ager and Luckhaus \cite{1992-TAMS-JaegerLuckhaus} consider a  parabolic--elliptic Keller--Segel system:
\begin{equation*}
\begin{cases}
\partial_t u=\Delta u-\operatorname{div}(u\nabla v),
    & x\in\Omega,\ t>0,\\[1mm]
0=\Delta v-\dfrac{M}{|\Omega|}+u,
    & x\in\Omega,\ t>0,\\[1mm]
\partial_\nu u=\partial_\nu v=0,
    & x\in\partial\Omega,\ t>0,\\[1mm]
u(\cdot,0)=u_0,
    & x\in\Omega,
\end{cases}
\end{equation*}
where 
 $\Omega\subset\mathbb R^2$ is a bounded domain with smooth boundary 
and $\nu$ is the outward unit normal to $\partial\Omega$ and $M=\int_\Omega u_0\dd x$.
J\"ager and Luckhaus \cite{1992-TAMS-JaegerLuckhaus} proved global existence for initial data of sufficiently small mass and finite-time blow-up for certain radially symmetric initial data.
Nagai \cite{1995-AMSA-Nagai} considered the parabolic--elliptic system with $0=\Delta v-v+u$ in a disk under radial symmetry, proving uniform boundedness for $M<8\pi$ and finite-time blow-up for sufficiently concentrated data of any mass $M>8\pi$.
Nagai et al. \cite{2001-RIMS-NagaiSenbaSuzuki} noted that solutions to the same system on $\Omega$ remain globally bounded when $M<4\pi$, without any symmetry assumption.
Nagai \cite{2001-JIA-Nagai} also proved finite-time blow-up when the initial data are sufficiently concentrated near an interior point and $M>8\pi$, or near a boundary point and $M>4\pi$ under suitable boundary assumptions.

For the parabolic--elliptic system with $0=\Delta v+u$ on $\mathbb R^2$, Blanchet et al. \cite{2006-EJDE-BlanchetDolbeaultPerthame} proved, under finite entropy and second moment assumptions, global existence of weak solutions for mass $M<8\pi$ and finite-time blow-up for mass $M>8\pi$.
Wei \cite{2018-JFA-Wei} later proved global existence of mild solutions for nonnegative initial data in $L^1(\mathbb R^2)$ with mass $M\le 8\pi$.

\subsection{Direct signal production on Riemannian surfaces}
Keller--Segel systems have also been studied on Riemannian surfaces.
Hillen and Painter \cite{2001-AAM-HillenPainter} proved global existence of classical solutions on closed Riemannian manifolds for a chemotaxis model preventing overcrowding.
In the parabolic--elliptic system with $0=\Delta_g v+u$ on the two-dimensional hyperbolic space,
Maheux and Pierfelice \cite{2020-SJMA-MaheuxPierfelice} proved global existence for mass $M<8\pi$ with finite weighted moment and finite-time blow-up for $M>8\pi$ with sufficiently small weighted moment. 
Nagy \cite{2023-PRSESA-Nagy} proved that, under suitable assumptions, the mass of positive stationary solutions on curved planes is necessarily $8\pi$.
For the parabolic--elliptic system with $0=\Delta_g v-v+u$, Li et al. \cite{LiLiYan-submitted} proved global boundedness for $M<8\pi$ on closed Riemannian surfaces and finite-time blow-up for sufficiently concentrated axially symmetric data with $M>8\pi$ on spheres and spheroids.

\subsection{Indirect signal production in planar domains}
In a bounded domain $\Omega\subset\mathbb R^2$, the system corresponding to \eqref{eq:system} arises from the mountain pine beetle model of Strohm et al. \cite{2013-BMB-StrohmTysonPowell} and takes the form
\begin{equation}\label{eq:planar-evolution}
\begin{cases}
\partial_t u=\nabla\cdot(\nabla u-u\nabla v),
&x\in\Omega,\ t>0,\\[1mm]
0=\Delta v-\overline w+w,
&x\in\Omega,\ t>0,\\[1mm]
\partial_t w=u-w,
&x\in\Omega,\ t>0,\\[1mm]
(\nabla u-u\nabla v)\cdot\nu=\partial_\nu v=0,
&x\in\partial\Omega,\ t>0,\\[1mm]
u(x,0)=u_0(x),\quad w(x,0)=w_0(x),
&x\in\Omega,
\end{cases}
\end{equation}
where $u$ and $w$ denote the densities of flying and nesting beetles, respectively, and $v$ is the pheromone produced by nesting beetles.
Set $A:=|\Omega|$ and $\overline w:=A^{-1}\int_\Omega w(x,t)\dd x$, and normalize $v$ to have zero spatial mean.
Write $m:=M/A$, where $M:=\int_\Omega u_0\dd x$ is the conserved cell mass.
In the unit disk, Tao and Winkler \cite{2017-JEMS-TaoWinkler} proved global classical solvability and, under radial symmetry, uniform boundedness for $M<8\pi$ and infinite-time blow-up for suitable initial data of any mass $M>8\pi$.
For the system \eqref{eq:planar-evolution} with $\partial_t v=\Delta v-v+w$, Lauren\c{c}ot \cite{2019-DCDSB-Laurencot} proved uniform boundedness for $M<4\pi$ and the existence of unbounded global solutions for $M>4\pi$ with $M\notin4\pi\mathbb N$ in smooth bounded domains in $\mathbb R^2$.
For the Cauchy problem in $\mathbb R^2$ with $0=\Delta v+w$, Xiang and Yang \cite{2025-JEE-XiangYang} proved global boundedness for $M<4\pi$, and, under radial symmetry, global boundedness for $M<8\pi$ and the existence of unbounded global solutions for $M>8\pi$.
For radial solutions of \eqref{eq:planar-evolution} in the unit disk, Mao and Li \cite{2024-DCDSS-MaoLi} proved Dirac-type chemotactic collapse with mass $M\ge32\pi$ at the origin in infinite time under suitable concentration assumptions.
When the second equation is replaced by $\partial_t v=\Delta v-v+w$ in a disk, Soga \cite{2025-JEE-Soga} showed that radially symmetric solutions blowing up in infinite time have a singularity like a Dirac delta function.

\subsection{Main results}
Here we study \eqref{eq:system} on closed connected Riemannian surfaces.
We prove global classical solvability, uniform boundedness for $M<8\pi$, and the existence of unbounded global solutions for every $M>8\pi$.
On spheres, we also establish uniform boundedness at $M=8\pi$ without assuming axial symmetry.
For notational convenience, we omit the domain $\Sigma$ from norm subscripts in what follows.
Our first result gives global existence.

\begin{proposition}\label{prop:global-existence}
Let $u_0\in C(\Sigma)$ and $w_0\in C^1(\Sigma)$ be nonnegative initial data.
Then \eqref{eq:system}--\eqref{eq:normalization} admits a unique global classical solution.
\end{proposition}

The following theorem gives uniform boundedness in the subcritical case.

\begin{theorem}\label{thm:subcritical}
Let $u_0\in C(\Sigma)$ and $w_0\in C^1(\Sigma)$ be nonnegative, and let $(u,v,w)$ be the corresponding global classical solution of \eqref{eq:system}--\eqref{eq:normalization}.
Assume that
\begin{equation*}
 M=\int_{\Sigma}u_0\dd\mu_g<8\pi.
\end{equation*}
Then there exists $C>0$, independent of $t$, such that
\begin{equation*}
 \norm{u(t)}_{L^\infty}
 +\norm{w(t)}_{L^\infty}
 +\norm{v(t)}_{W^{1,\infty}}
 \le C
 \quad\text{for all }t\ge0.
\end{equation*}
\end{theorem}

The next theorem establishes the existence of solutions exhibiting infinite-time blow-up for every $M>8\pi$.

\begin{theorem}\label{thm:supercritical}
Let $M>8\pi$.
There exist strictly positive $u_0,w_0\in C^\infty(\Sigma)$ with $\int_\Sigma u_0\dd\mu_g=M$ such that the corresponding global classical solution of \eqref{eq:system}--\eqref{eq:normalization} satisfies
\begin{equation*}
 \limsup_{t\to\infty}\norm{u(t)}_{L^\infty}=\infty.
\end{equation*}
\end{theorem}
The same argument can be used to deal with the supercritical case of \eqref{eq:planar-evolution} by extending $\Omega$ to a compact connected Riemann surface.
\begin{theorem}\label{thm:planar-infinite-time-blow-up}
For every $M>4\pi$, there exist $u_0,w_0\in C^\infty(\overline\Omega)$ with
$$\int_\Omega u_0\dd x=M$$ such that the corresponding solution $(u, v, w)$ of \eqref{eq:planar-evolution} satisfies
\[
\limsup_{t\to\infty}\norm{u(\cdot,t)}_{L^\infty(\Omega)}=\infty.
\]
\end{theorem}

Finally, we consider the critical mass $M=8\pi$. 
On spheres equipped with the canonical spherical metric, the solution is globally bounded.
We also
construct a rotationally invariant solution that blows up in infinite time on
the Riemann sphere equipped with a metric
defined in action--angle coordinates as in \cite{2003-EnseignMath-Hwang}.

\begin{theorem}\label{thm:critical-sphere}
Let $R>0$, and let $\Sigma=\mathbb S_R:=\{x\in\mathbb R^3:|x|=R\}$ be endowed with the canonical spherical metric \cite[pp.~15--16]{2018-Lee}.
Let $u_0\in C(\Sigma)$ and $w_0\in C^1(\Sigma)$ be nonnegative, and assume that
\[
 \int_\Sigma u_0\dd\mu_g=8\pi.
\]
Then there exists $C>0$, independent of $t$, such that the solution of \eqref{eq:system}--\eqref{eq:normalization} satisfies
\begin{equation*}
 \norm{u(t)}_{L^\infty}
 +\norm{w(t)}_{L^\infty}
 +\norm{v(t)}_{W^{1,\infty}}
 \le C
 \quad\text{for all }t\ge0.
\end{equation*}
\end{theorem}

\begin{theorem}\label{thm:critical-grow-up}
Let $\Sigma:=\mathbb C\cup\{\infty\}$, with the action given by $z\mapsto e^{i\beta}z$ for $\beta\in\mathbb R$.
For $z\in\mathbb C\setminus\{0\}$, set $\theta:=\arg z$ and define $s\in(0,1)$ by
\begin{equation}\label{eq:critical-grow-up-coordinate}
 \log|z|=\int_{1/2}^{s}\frac{\dd\sigma}{2D(\sigma)},\quad
 D(s):=s(1-s)\left(1-\tfrac12s(1-s)\right).
\end{equation}
Equip $\Sigma$ with the smooth extension of the metric
\[
 g:=D(s)^{-1}\mathrm{d}s^2+4D(s)\,\mathrm{d}\theta^2.
\]
Then there exist strictly positive, rotation-invariant initial data $u_0=w_0\in C^\infty(\Sigma)$ such that
\[
 \int_\Sigma u_0\dd\mu_g=\int_\Sigma w_0\dd\mu_g=8\pi
\]
and the corresponding global classical solution of \eqref{eq:system}--\eqref{eq:normalization} satisfies
\[
 \norm{u(t)}_{L^\infty}\ge 1000+\frac{t}{2000}
 \qquad\text{for all }t\ge0.
\]
\end{theorem}

We briefly outline the main arguments used to prove these results.
Proposition~\ref{prop:global-existence} follows from estimates as in \cite[Lemmas~3.1 and 3.2]{2017-JEMS-TaoWinkler}, combined with heat semigroup estimates.
The key step in the proof of Theorem~\ref{thm:subcritical} is to control the entropy using a Lyapunov functional and the Moser--Trudinger inequality \cite{2007-TAMS-Yang}.
The proof of Theorem~\ref{thm:supercritical} relies on a uniform lower bound for the energy of stationary solutions, obtained from the pointwise estimates in \cite[Theorem~0.2]{1999-CMP-Li}.
Both critical-mass results rely on one-dimensional comparison.
For Theorem~\ref{thm:critical-sphere}, we use spherical isoperimetry \cite{1978-BAMS-Osserman} and construct a supersolution as in \cite[Lemma~4.5]{2025-JEE-XiangYang}.
For Theorem~\ref{thm:critical-grow-up}, we use the comparison principle and adapt the subsolution profile from \cite{2017-JEMS-TaoWinkler}.

The paper is organized as follows.
Section~\ref{sec:semigroup} contains the heat semigroup estimates.
Local and global existence are proved in Sections~\ref{sec:local-existence} and~\ref{sec:global-existence}, respectively.
Sections~\ref{sec:subcritical} and~\ref{sec:supercritical} establish uniform boundedness for subcritical mass and construct unbounded solutions for supercritical mass, respectively.
Section~\ref{sec:planar-neumann} proves the planar Neumann result.
At critical mass, Sections~\ref{sec:critical-sphere} and~\ref{sec:critical-grow-up} prove uniform boundedness on spheres and construct a solution that blows up in infinite time on a closed surface, respectively.

\section{Heat semigroup estimates}
\label{sec:semigroup}
In this section, we establish heat semigroup estimates on $(\Sigma,g)$.
For $p\in(1,\infty)$ and $s\in\mathbb R$, we denote the Sobolev spaces on $\Sigma$ by $W^{s,p}(\Sigma)$ (see \cite[Section~13.6]{2011-Taylor}).
We first establish Gaussian upper bounds for the heat kernel and its gradient.

\begin{lemma}\label{lem:heat-kernel-bounds}
Let $H(x,y,\tau)$ be the heat kernel of $(\Sigma,g)$.
There exist $A_1,A_2,A_3,A_4>0$ such that
\begin{equation}\label{eq:heat-kernel-gaussian-bound}
 H(x,y,\tau)\le A_1\tau^{-1}
 \exp\left(-A_2\frac{d_g(x,y)^2}{\tau}\right)
 \quad \text{for all } x,y\in\Sigma,\ \tau\in(0,2],
\end{equation}
and
\begin{equation}\label{eq:heat-kernel-gradient}
 \left|\nabla_{g,y}H(x,y,\tau)\right|_g
 \le A_3\tau^{-\frac{3}{2}}
 \exp\left(-A_4\frac{d_g(x,y)^2}{\tau}\right)
 \quad \text{for all } x,y\in\Sigma,\ \tau\in(0,1].
\end{equation}
\end{lemma}

\begin{proof}
By \cite[Theorem~15.14]{2009-Grigoryan} and the compactness of $\Sigma$, there are $C_1,C_2>0$ such that \eqref{eq:heat-kernel-gaussian-bound} holds with $A_1=C_1$ and $A_2=C_2$.
To prove \eqref{eq:heat-kernel-gradient}, we choose $a\in(0,\frac{1}{2})$ smaller than the injectivity radius of $(\Sigma,g)$, set $\ell=a\sqrt{\tau}$, and consider the cylinder
\begin{equation*}
 Q_\ell=B_g(y,\ell)\times(\tau-\ell^2,\tau+\ell^2).
\end{equation*}
In normal coordinates centered at $y$, introduce the rescaled function
\begin{equation*}
 \mathcal H(\xi,\theta)
 :=H\bigl(x,\exp_y(\ell\xi),\tau+\ell^2\theta\bigr)
 \quad \text{for all } (\xi,\theta)\in B_1(0)\times(-1,1).
\end{equation*}
Since $H(x,z,s)$ satisfies $(\partial_s-\Delta_{g,z})H(x,z,s)=0$ in $Q_\ell$, the interior estimates in \cite[Theorem~6.21]{2009-Grigoryan} on smaller cylinders containing $(0,0)$ and the Sobolev embedding in \cite[Theorem~6.1]{2009-Grigoryan} give $C_3,C_4>0$ such that
\begin{equation*}
 |\nabla_\xi\mathcal H(0,0)|
 \le C_3\norm{\mathcal H}_{L^2(B_1(0)\times(-1,1))}
 \le C_4\sup_{(z,s)\in Q_\ell}H(x,z,s)
 \quad\text{for all }x,y\in\Sigma,\ 0<\tau\le1,
\end{equation*}
which yields
\begin{equation}\label{eq:heat-kernel-local-gradient}
 |\nabla_{g,y}H(x,y,\tau)|_g
 \le C_4\ell^{-1}\sup_{(z,s)\in Q_\ell}H(x,z,s)
 \quad\text{for all }x,y\in\Sigma,\ 0<\tau\le1.
\end{equation}
For $(z,s)\in Q_\ell$, we have $(1-a^2)\tau<s<(1+a^2)\tau$ and $d_g(y,z)<\ell$.
Consequently, for all $(z,s)\in Q_\ell$,
\begin{equation*}
 \frac{d_g(x,z)^2}{s}
 \ge\frac{d_g(x,y)^2}{2s}-\frac{\ell^2}{s}
 \ge\frac{d_g(x,y)^2}{2(1+a^2)\tau}-\frac{a^2}{1-a^2}.
\end{equation*}
Since $s\in(0,(1+a^2)\tau)\subset(0,2]$, the Gaussian bound \eqref{eq:heat-kernel-gaussian-bound} gives
\[
 \sup_{(z,s)\in Q_\ell}H(x,z,s)
 \le C_5\tau^{-1}
 \exp\left(-C_6\frac{d_g(x,y)^2}{\tau}\right)
 \quad\text{for all }x,y\in\Sigma,\ 0<\tau\le1,
\]
with $C_5:=C_1(1-a^2)^{-1}\exp(\frac{C_2a^2}{1-a^2})$ and $C_6:=\frac{C_2}{2(1+a^2)}$.
Combining this with \eqref{eq:heat-kernel-local-gradient} and $\ell=a\sqrt{\tau}$ gives \eqref{eq:heat-kernel-gradient} with $A_3=C_4C_5a^{-1}$ and $A_4=C_6$.
\end{proof}

We next derive $L^q$--$L^p$ smoothing estimates for the heat semigroup.

\begin{lemma}\label{lem:semigroup-smoothing}
Let $1\le q\le p\le\infty$.
There exists a constant $C=C(\Sigma,g,p,q)>0$ such that, for all $\tau\in(0,1]$, $f\in L^q(\Sigma)$, and $F\in L^q(T\Sigma)$,
\begin{align*}
 \norm{e^{\tau\Delta_g}f}_{L^p}
 &\le C\tau^{-\left(\frac{1}{q}-\frac{1}{p}\right)}\norm{f}_{L^q},\\
 \norm{e^{\tau\Delta_g}\operatorname{div}_gF}_{L^p}
 &\le C\tau^{-\frac{1}{2}-\left(\frac{1}{q}-\frac{1}{p}\right)}
 \norm{F}_{L^q},
\end{align*}
where $\operatorname{div}_gF$ is understood in the sense of distributions.
\end{lemma}

\begin{proof}
By Lemma~\ref{lem:heat-kernel-bounds}, choose $C_1>0$ sufficiently large and set
\[
 K_\tau(x,y):=\tau^{-1}
 \exp\left(-\frac{d_g(x,y)^2}{C_1\tau}\right),
 \qquad
 T_\tau h(x):=\int_\Sigma K_\tau(x,y)h(y)\dd\mu_g(y).
\]
The heat kernel representation \cite[Theorem~7.13]{2009-Grigoryan} and its extension to distributions \cite[(5.4)]{2015-TAMS-KerkyacharianPetrushev} give
\[
 \bigl(e^{\tau\Delta_g}\operatorname{div}_gF\bigr)(x)
 =-\int_\Sigma
 \langle\nabla_{g,y}H(x,y,\tau),F(y)\rangle_g\dd\mu_g(y).
\]
Hence \eqref{eq:heat-kernel-gaussian-bound} and \eqref{eq:heat-kernel-gradient} imply
\begin{align}
 |e^{\tau\Delta_g}f(x)|
 &\le C_1T_\tau(|f|)(x),\notag\\
 \left|\bigl(e^{\tau\Delta_g}\operatorname{div}_gF\bigr)(x)\right|
 &\le C_1\tau^{-\frac12}T_\tau(|F|_g)(x)
 \quad\text{for all }x\in\Sigma,\ 0<\tau\le1.
 \label{eq:semigroup-pointwise-bounds}
\end{align}
It therefore suffices to estimate $T_\tau$.
For $1\le r<\infty$, geodesic polar coordinates and the substitution $\rho=\sqrt{\tau}\,s$ give $C_2,C_3>0$ such that
\begin{align*}
 \int_\Sigma K_\tau(x,y)^r\dd\mu_g(y)
 &\le C_2\tau^{-r}\int_0^{\operatorname{diam}(\Sigma)}
       e^{-r\rho^2/(C_1\tau)}\rho\dd\rho\\
 &=C_2\tau^{-r+1}\int_0^{\operatorname{diam}(\Sigma)/\sqrt{\tau}}
       e^{-rs^2/C_1}s\dd s\\
 &\le C_2\tau^{-r+1}\int_0^\infty
       e^{-rs^2/C_1}s\dd s\\
 &=\frac{C_1C_2}{2r}\tau^{-r+1}
 =:C_3\tau^{-r+1}
 \quad\text{for all }x\in\Sigma,\ 0<\tau\le1.
\end{align*}
Thus, with $C_4:=C_3^{1/r}$,
\[
 \sup_{x\in\Sigma}\norm{K_\tau(x,\cdot)}_{L^r}
 =\sup_{y\in\Sigma}\norm{K_\tau(\cdot,y)}_{L^r}
 \le C_4\tau^{-1+\frac1r}
 \quad\text{for all }0<\tau\le1.
\]
Let $r'$ be the conjugate exponent of $r$.
For $h\in L^1(\Sigma)$, Minkowski's integral inequality \cite[Theorem~6.19(b)]{1999-Folland} gives
\begin{align*}
 \norm{T_\tau h}_{L^r}
 &\le\int_\Sigma\norm{K_\tau(\cdot,y)}_{L^r}|h(y)|\dd\mu_g(y)\\
 &\le C_4\tau^{-1+\frac1r}\norm{h}_{L^1}
 \quad\text{for all }0<\tau\le1.
\end{align*}
For $h\in L^{r'}(\Sigma)$, H\"older's inequality yields
\begin{align*}
 |T_\tau h(x)|
 &\le\norm{K_\tau(x,\cdot)}_{L^r}\norm{h}_{L^{r'}}\\
 &\le C_4\tau^{-1+\frac1r}\norm{h}_{L^{r'}}
 \quad\text{for all }x\in\Sigma,\ 0<\tau\le1.
\end{align*}
Taking the supremum over $x$ in the latter estimate, we conclude that
\[
 \norm{T_\tau}_{L^1\to L^r},\quad
 \norm{T_\tau}_{L^{r'}\to L^\infty}
 \le C_4\tau^{-1+\frac1r}
 \quad\text{for all }0<\tau\le1.
\]
For $1\le q\le p<\infty$, choose $r$ and $\theta$ by
\[
 \frac1r=1-\frac1q+\frac1p,
 \qquad \theta=1-\frac rp\in[0,1).
\]
The Riesz--Thorin interpolation theorem \cite[Theorem~6.27]{1999-Folland}, or the $L^1\to L^r$ bound when $\theta=0$, gives
\[
 \norm{T_\tau h}_{L^p}
 \le C_4\tau^{-\left(\frac1q-\frac1p\right)}\norm{h}_{L^q}
 \quad\text{for all }0<\tau\le1.
\]
For $p=\infty$ and $q>1$, the same estimate follows from the $L^{r'}\to L^\infty$ bound with $r=q'$.
The remaining case $p=\infty$, $q=1$ follows directly from $K_\tau\le\tau^{-1}$.
Combining these estimates with \eqref{eq:semigroup-pointwise-bounds} completes the proof.
\end{proof}

\section{Local existence}
\label{sec:local-existence}
We first establish local existence and uniqueness of classical solutions, together with a blow-up criterion.
\begin{proposition}\label{prop:local-existence}
For nonnegative initial data $u_0\in C(\Sigma)$ and $w_0\in C^1(\Sigma)$, there exists $T_{\max}\in(0,\infty]$ such that \eqref{eq:system}--\eqref{eq:normalization} admits a unique classical solution $(u,v,w)$ satisfying
\begin{align}
 u&\in C(\Sigma\times[0,T_{\max}))
 \cap C^{2,1}(\Sigma\times(0,T_{\max})),\notag\\
 v&\in C([0,T_{\max});C^2(\Sigma)),\notag\\
 w&\in C^1([0,T_{\max});C(\Sigma)).
 \label{eq:local-solution-regularity}
\end{align}
Moreover, $u,w\ge0$, and
\begin{equation}\label{eq:local-mass-k}
 \overline w=\frac{e^{-t}}{\mu_g(\Sigma)}
       \int_\Sigma w_0\dd\mu_g
       +\frac{M}{\mu_g(\Sigma)}(1-e^{-t})
 \quad \text{for all } t\in[0,T_{\max}).
\end{equation}
If $T_{\max}<\infty$, then
\begin{equation}\label{eq:local-blow-up-alternative}
 \limsup_{t\nearrow T_{\max}}\norm{u(t)}_{L^\infty}=\infty.
\end{equation}
\end{proposition}

\begin{proof}
\medskip\noindent\textit{Step 1. Existence.}
Let $T\in(0,1]$ and define
\[
 X_T:=C([0,T];C(\Sigma)),\quad
 \norm{f}_{X_T}:=\sup_{t\in[0,T]}\norm{f(t)}_{L^\infty},
 \quad R:=2\norm{u_0}_{L^\infty}+1.
\]
Consider the closed ball $S_T:=\{f\in X_T:\norm{f}_{X_T}\le R\}$.
For $f\in S_T$, define
\begin{equation}\label{eq:local-w-v}
 w_f(t):=e^{-t}w_0+\int_0^t e^{-(t-s)}f(s)\dd s.
\end{equation}
For each $t\in[0,T]$, let $v_f(t)$ be the unique solution of
\[
 -\Delta_gv_f(t)
 =w_f(t)-\frac{1}{\mu_g(\Sigma)}\int_\Sigma w_f(t)\dd\mu_g,
 \quad \int_\Sigma v_f(t)\dd\mu_g=0.
\]
The definition of $w_f$ in \eqref{eq:local-w-v} implies
\[
 \norm{w_f(t)}_{L^\infty}
 \le\norm{w_0}_{L^\infty}+\norm{f}_{X_T}
 \le\norm{w_0}_{L^\infty}+R
 \quad \text{for all } t\in(0,T].
\]
Fix $p>2$.
The elliptic estimate \cite[Theorem A.3]{2019-Lee} and Sobolev embedding \cite[Theorem~2.7]{1999-Hebey} yield constants $C_1,C_2>0$ such that
\begin{align}
 \norm{v_f(t)}_{W^{2,p}}
 +\norm{\nabla_gv_f(t)}_{L^\infty}
 &\le C_1\left\|w_f(t)
 -\frac{1}{\mu_g(\Sigma)}\int_\Sigma w_f(t)\dd\mu_g\right\|_{L^p}\notag\\
 &\le C_2\norm{w_f(t)}_{L^p}\notag\\
 &\le C_3(\norm{w_0}_{L^\infty}+R)
 \quad \text{for all } t\in(0,T],\label{eq:local-elliptic}
\end{align}
where $C_3:=C_2\mu_g(\Sigma)^{\frac{1}{p}}$.
With $B:=\norm{w_0}_{L^\infty}+R$, these estimates imply
\begin{equation}\label{eq:local-candidate-bounds}
 \sup_{t\in[0,T]}\norm{w_f(t)}_{L^\infty}\le B,
 \quad
 \sup_{t\in[0,T]}\norm{\nabla_gv_f(t)}_{L^\infty}\le C_3B.
\end{equation}
For $f_1,f_2\in S_T$, \eqref{eq:local-w-v} gives
\begin{equation}\label{eq:local-w-difference}
 \sup_{t\in[0,T]}\norm{w_{f_1}(t)-w_{f_2}(t)}_{L^\infty}
 \le T\norm{f_1-f_2}_{X_T}.
\end{equation}
Since $v_{f_1}-v_{f_2}$ has zero mean and satisfies
\[
 -\Delta_g(v_{f_1}-v_{f_2})
 =w_{f_1}-w_{f_2}
 -\frac{1}{\mu_g(\Sigma)}\int_\Sigma(w_{f_1}-w_{f_2})\dd\mu_g,
\]
the elliptic estimate \cite[Theorem~A.3]{2019-Lee} and \eqref{eq:local-w-difference} yield
\begin{equation}\label{eq:local-candidate-difference}
\sup_{t\in[0,T]}\norm{\nabla_g(v_{f_1}-v_{f_2})(t)}_{L^\infty}
 \le C_3T\norm{f_1-f_2}_{X_T}.
\end{equation}
Define
\begin{equation}\label{eq:local-fixed-point}
 (\Phi f)(t):=e^{t\Delta_g}u_0
 -\int_0^t e^{(t-s)\Delta_g}
       \operatorname{div}_g(f(s)\nabla_gv_f(s))\dd s.
\end{equation}
Let $\delta:=\frac{1}{2}-\frac{1}{p}>0$.
Lemma~\ref{lem:semigroup-smoothing} gives a constant $C_4>0$ such that
\begin{equation}\label{eq:local-div-semigroup}
 \norm{e^{\tau\Delta_g}\operatorname{div}_gF}_{L^\infty}
 \le C_4\tau^{-1+\delta}\norm{F}_{L^p}\quad \text{for all } \tau\in(0,1].
\end{equation}
Let $C_5:=\mu_g(\Sigma)^{\frac{1}{p}}C_3$ and $C_6:=C_4C_5$.
The contractivity of the heat semigroup on $C(\Sigma)$ \cite[Exercise~7.33]{2009-Grigoryan}, together with \eqref{eq:local-candidate-bounds}, \eqref{eq:local-fixed-point} and \eqref{eq:local-div-semigroup}, gives
\begin{align*}
 \norm{(\Phi f)(t)}_{L^\infty}
 &\le \norm{e^{t\Delta_g}u_0}_{L^\infty}
 +\int_0^t\left\|e^{(t-s)\Delta_g}
    \operatorname{div}_g(f(s)\nabla_gv_f(s))\right\|_{L^\infty}\dd s\\
 &\le \norm{u_0}_{L^\infty}
 +C_4\int_0^t(t-s)^{-1+\delta}
    \norm{f(s)\nabla_gv_f(s)}_{L^p}\dd s\\
 &\le \norm{u_0}_{L^\infty}
 +C_4\mu_g(\Sigma)^{\frac{1}{p}}\int_0^t(t-s)^{-1+\delta}
    \norm{f(s)}_{L^\infty}\norm{\nabla_gv_f(s)}_{L^\infty}\dd s\\
 &\le \norm{u_0}_{L^\infty}
 +C_6RB\int_0^t(t-s)^{-1+\delta}\dd s\\
 &=\norm{u_0}_{L^\infty}+C_6RB\frac{t^\delta}{\delta}
 \quad \text{for all } t\in(0,T].
\end{align*}
The above estimates and the strong continuity of the heat semigroup imply that $\Phi f\in X_T$ for every $f\in S_T$, with
\begin{equation}\label{eq:local-self-map}
 \norm{\Phi f}_{X_T}
 \le\norm{u_0}_{L^\infty}+C_6RB\frac{T^\delta}{\delta}.
\end{equation}
Furthermore, by \eqref{eq:local-candidate-bounds}, \eqref{eq:local-candidate-difference}, \eqref{eq:local-fixed-point}, \eqref{eq:local-div-semigroup}, and $\norm{f_2}_{X_T}\le R$, we obtain
\begin{align}
 \norm{(\Phi f_1)(t)-(\Phi f_2)(t)}_{L^\infty}
 &\le \int_0^t\left\|e^{(t-s)\Delta_g}\operatorname{div}_g
   \bigl(f_1(s)\nabla_gv_{f_1}(s)
         -f_2(s)\nabla_gv_{f_2}(s)\bigr)\right\|_{L^\infty}\dd s\notag\\
 &\le C_4\int_0^t(t-s)^{-1+\delta}
   \Bigl(\norm{(f_1-f_2)(s)\nabla_gv_{f_1}(s)}_{L^p}\notag\\
 &\hspace{16mm}
       +\norm{f_2(s)\nabla_g(v_{f_1}-v_{f_2})(s)}_{L^p}\Bigr)\dd s\notag\\
 &\le C_4\mu_g(\Sigma)^{\frac1p}\int_0^t(t-s)^{-1+\delta}
   \Bigl(\norm{(f_1-f_2)(s)}_{L^\infty}
          \norm{\nabla_gv_{f_1}(s)}_{L^\infty}\notag\\
 &\hspace{16mm}
       +\norm{f_2(s)}_{L^\infty}
          \norm{\nabla_g(v_{f_1}-v_{f_2})(s)}_{L^\infty}\Bigr)\dd s\notag\\
 &\le C_6(B+RT)\norm{f_1-f_2}_{X_T}
       \int_0^t(t-s)^{-1+\delta}\dd s\notag\\
 &=C_6(B+RT)\frac{t^\delta}{\delta}\norm{f_1-f_2}_{X_T}\notag\\
&\le C_6(B+RT)\frac{T^\delta}{\delta}
       \norm{f_1-f_2}_{X_T}  \quad \text{for all } t\in(0,T].
 \label{eq:local-contraction}
\end{align}
Choose $T\in(0,1]$ sufficiently small that $\frac{C_6RBT^\delta}{\delta}\le R-\norm{u_0}_{L^\infty}$ and $\frac{C_6(B+RT)T^\delta}{\delta}\le\frac{1}{2}$.
By \eqref{eq:local-self-map} and \eqref{eq:local-contraction}, the Banach fixed point theorem yields a unique fixed point $u\in S_T$.
Together with $w=w_u$, $v=v_u$, and $\overline w$ defined by \eqref{eq:kdef}, this yields a mild solution of \eqref{eq:system}.

\medskip\noindent\textit{Step 2. Regularity.}
The formula \eqref{eq:local-w-v} gives $w\in C^1([0,T];C(\Sigma))$, while \eqref{eq:local-elliptic} yields $v\in C([0,T];W^{2,p}(\Sigma))$ for all $p>2$.
Fix $\alpha\in \left(0,1-\frac{2}{p}\right)$ and choose $\varepsilon>0$ sufficiently small that $\beta:=\frac{1+\alpha}{2}+\frac{1}{p}+\frac{\varepsilon}{2}<1$.
For $F\in L^p(T\Sigma)$ and $\varphi\in W^{1,p'}(\Sigma)$, where $p'=\frac{p}{p-1}$, H\"older's inequality yields
\[
 |\langle\operatorname{div}_gF,\varphi\rangle|
 \le\norm{F}_{L^p}\norm{\nabla_g\varphi}_{L^{p'}}
 \le\norm{F}_{L^p}\norm{\varphi}_{W^{1,p'}}.
\]
Thus, by the dual characterization of $W^{-1,p}(\Sigma)$ \cite[Theorem~100]{2022-M-BehzadanHolst}, $\operatorname{div}_g:L^p(T\Sigma)\to W^{-1,p}(\Sigma)$ is bounded.
By the heat semigroup estimates \cite[(15.1.13) and (15.1.15)]{2011-Taylor} and the Sobolev embedding \cite[pp.~43--44]{2011-Taylor}, we can combine \eqref{eq:local-candidate-bounds}, \eqref{eq:local-fixed-point}, and $\norm{u}_{X_T}\le R$ to find $C_7\ge1$ and $C_8>0$ such that
\begin{align}
 \norm{u(t)}_{C^{\alpha}}
 &\le\norm{e^{t\Delta_g}u_0}_{C^{\alpha}}
   +\int_0^t\norm{e^{(t-s)\Delta_g}
       \operatorname{div}_g(u(s)\nabla_gv(s))}_{C^{\alpha}}\dd s\notag\\
 &\le C_7t^{-\frac{\alpha}{2}}\norm{u_0}_{L^\infty}\notag\\
 &\quad+C_7\int_0^t
       \norm{e^{(t-s)\Delta_g}\operatorname{div}_g(u(s)\nabla_gv(s))}
       _{W^{\alpha+\frac2p+\varepsilon,p}}\dd s\notag\\
 &\le C_7t^{-\frac{\alpha}{2}}\norm{u_0}_{L^\infty}
   +C_7^2\int_0^t(t-s)^{-\beta}
       \norm{\operatorname{div}_g(u(s)\nabla_gv(s))}_{W^{-1,p}}\dd s\notag\\
 &\le C_7t^{-\frac{\alpha}{2}}\norm{u_0}_{L^\infty}
   +C_7^3\int_0^t(t-s)^{-\beta}
       \norm{u(s)\nabla_gv(s)}_{L^p}\dd s\notag\\
 &\le C_7t^{-\frac{\alpha}{2}}\norm{u_0}_{L^\infty}
   +C_7^3C_5RB\int_0^t(t-s)^{-\beta}\dd s\notag\\
 &=C_7t^{-\frac{\alpha}{2}}\norm{u_0}_{L^\infty}
   +\frac{C_7^3C_5RB}{1-\beta}t^{1-\beta}\notag\\
 &\le C_7t^{-\frac{\alpha}{2}}\norm{u_0}_{L^\infty}
   +\frac{C_7^3C_5RB}{1-\beta}\notag\\
 &\le C_8\left(1+t^{-\frac{\alpha}{2}}\right)
 \quad\text{for all }t\in(0,T].
 \label{eq:local-u-holder-bound}
\end{align}
Since $w_0\in C^1(\Sigma)$, \eqref{eq:local-w-v} implies that $w\in C([0,T];C^{\alpha}(\Sigma))$.
Elliptic Schauder regularity (see \cite[Theorem A.6]{2019-Lee}) yields $v\in C([0,T];C^{2+\alpha}(\Sigma))$.
Since $\partial_t w=u-w$ and $u,w$ are bounded, there exists $C_{9}>0$ such that
\begin{align*}
 \norm{w(t)-w(s)}_C
 &\le \int_s^t\norm{u(r)-w(r)}_C\dd r\\
 &\le (t-s)\sup_{r\in[0,T]}
 \bigl(\norm{u(r)}_C
       +\norm{w(r)}_C\bigr)\\
 &\le C_{9}|t-s|
 \quad\text{for all }0\le s<t\le T.
\end{align*}
Moreover,
\[
 \sup_{t\in[0,T]}\norm{w(t)}_{C^{\alpha}}\le C_{9}.
\]
Choose $\theta\in\left(0,\frac{2\alpha}{2+\alpha}\right)$.
Interpolation gives
\begin{align*}
 \norm{w(t)-w(s)}_{C^{\theta}}
 &\le C_{10}\norm{w(t)-w(s)}_C^{1-\frac{\theta}{\alpha}}
       \norm{w(t)-w(s)}_{C^{\alpha}}^{\frac{\theta}{\alpha}}\\
 &\le C_{11}|t-s|^{1-\frac{\theta}{\alpha}}
 \le C_{11}|t-s|^{\frac{\theta}{2}}
 \quad\text{for all }0\le s<t\le T.
\end{align*}
Thus $w\in C^{\theta,\frac{\theta}{2}}(\Sigma\times[0,T])$.
Elliptic Schauder estimates \cite[Theorem~A.6]{2019-Lee} then imply that $\nabla_gv\in C^{\theta,\frac{\theta}{2}}(\Sigma\times[0,T])$.
Fix $\eta\in(0,T)$.
Using the semigroup property, the Gaussian bound \eqref{eq:heat-kernel-gaussian-bound}, \eqref{eq:local-fixed-point}, \eqref{eq:local-div-semigroup}, and the uniform $C^{\alpha}(\Sigma)$ bound in \eqref{eq:local-u-holder-bound}, we find constants $C_{12},C_{13},C_{14}>0$ such that
\begin{align*}
 \norm{u(t+h)-u(t)}_C
 &\le \norm{(e^{h\Delta_g}-I)u(t)}_C
       +\int_t^{t+h}\left\|e^{(t+h-s)\Delta_g}
       \operatorname{div}_g(u(s)\nabla_gv(s))\right\|_C\dd s\\
 &\le \norm{(e^{h\Delta_g}-I)u(t)}_C
       +C_4\int_t^{t+h}(t+h-s)^{-1+\delta}
       \norm{u(s)\nabla_gv(s)}_{L^p}\dd s\\
 &\le \norm{(e^{h\Delta_g}-I)u(t)}_C
       +C_6RB\int_t^{t+h}(t+h-s)^{-1+\delta}\dd s\\
 &\le \norm{u(t)}_{C^{\alpha}}\sup_{x\in\Sigma}
       \int_\Sigma H(x,y,h)d_g(x,y)^\alpha\dd\mu_g(y)
       +\frac{C_6RB}{\delta}h^\delta\\
 &\le C_{12}h^{-1}\sup_{x\in\Sigma}
       \int_\Sigma e^{-\frac{C_{13}d_g(x,y)^2}{h}}
       d_g(x,y)^\alpha\dd\mu_g(y)
       +\frac{C_6RB}{\delta}h^\delta\\
 &\le C_{14}h^{-1}\int_0^{\operatorname{diam}(\Sigma)}
       e^{-\frac{C_{13}\rho^2}{h}}\rho^{\alpha+1}\dd\rho
       +\frac{C_6RB}{\delta}h^\delta\\
 &\le C_{14}h^{\frac{\alpha}{2}}\int_0^\infty
       e^{-C_{13}r^2}r^{\alpha+1}\dd r
       +\frac{C_6RB}{\delta}h^\delta\\
 &=C_{15}h^{\frac{\alpha}{2}}+\frac{C_6RB}{\delta}h^\delta
 \le C_{16}h^{\frac{\alpha}{2}}
\end{align*}
for all $\eta/2\le t<T$ and $0<h\le T-t$.
Here we set $C_{15}:=C_{14}\int_0^\infty e^{-C_{13}r^2}r^{\alpha+1}\dd r$ and $C_{16}:=C_{15}+C_6RB/\delta$.
Hence, $u\in C^{\alpha,\frac{\alpha}{2}}(\Sigma\times\left[\frac{\eta}{2},T\right])$.
Choose a smooth time cutoff $\chi$ vanishing on $\left[0,\frac{\eta}{2}\right]$ and equal to $1$ on $[\eta,T]$.
Then $z:=\chi u$ satisfies
\begin{equation*}
\partial_t z-\Delta_g z+\left\langle\nabla_g v, \nabla_g z\right\rangle_g-(w-\overline w) z=\chi^{\prime} u,
\end{equation*}
in the sense of distributions.
By \cite[Theorem~2.4]{2026-arXiv-Huang}, this equation admits a classical solution $\widetilde z\in C^{2+\theta,1+\frac{\theta}{2}}(\Sigma\times[\frac{\eta}{2},T])$ with $\widetilde z\left(\frac{\eta}{2}\right)=0$.
Both $z$ and $\widetilde z$ satisfy the corresponding variation-of-constants formula.
Hence, by \eqref{eq:local-candidate-bounds} and \eqref{eq:local-div-semigroup},
\[
 \sup_{t\in[a,a+\ell]}\norm{z(t)-\widetilde z(t)}_C
 \le \frac{C_{17}\ell^\delta}{\delta}
       \sup_{t\in[a,a+\ell]}\norm{z(t)-\widetilde z(t)}_C
\]
whenever $[a,a+\ell]\subset[\eta/2,T]$ and $z(a)=\widetilde z(a)$, where $C_{17}:=C_4\mu_g(\Sigma)^{1/p}C_3B$.
Choosing $\ell>0$ with $C_{17}\ell^\delta<\delta$ and iterating from $z(\eta/2)=\widetilde z(\eta/2)=0$, we obtain $z=\widetilde z$ on $[\eta/2,T]$.
Since $\chi=1$ on $[\eta,T]$, it follows that $u\in C^{2,1}(\Sigma\times[\eta,T])$.
Since $\eta>0$ is arbitrary, we obtain $u\in C^{2,1}(\Sigma\times(0,T])$.
This establishes \eqref{eq:local-solution-regularity}.

\medskip\noindent\textit{Step 3. Uniqueness and nonnegativity.}
Let $(u_i,v_i,w_i)$, $i=1,2$, be two classical solutions with the same initial values.
Fix $T_*>0$ such that both solutions exist on $[0,T_*]$, and let
\[
 R_*:=1+\max_{i=1,2}\sup_{t\in[0,T_*]}\norm{u_i(t)}_{L^\infty},
 \quad
 B_*:=R_*+\max_{i=1,2}\sup_{t\in[0,T_*]}\norm{w_i(t)}_{L^\infty}.
\]
Whenever $(u_1,w_1)(a)=(u_2,w_2)(a)$, \eqref{eq:local-contraction} implies that for $\ell\in(0,1)$ and $[a,a+\ell]\subset[0,T_*]$,
\[
 \sup_{t\in[a,a+\ell]}\norm{u_1(t)-u_2(t)}_{L^\infty}
 \le C_6(B_*+R_*\ell)\frac{\ell^\delta}{\delta}\sup_{t\in[a,a+\ell]}\norm{u_1(t)-u_2(t)}_{L^\infty}.
\]
Choose $\ell_0\in(0,1)$ such that $\frac{C_6(B_*+R_*\ell_0)\ell_0^\delta}{\delta}<1$.
Then, for $\ell\in(0,\ell_0]$, we obtain $u_1=u_2$ on $[a,a+\ell]$, which yields $(u_1,v_1,w_1)=(u_2,v_2,w_2)$ for $t\in [a,a+\ell]$.
Starting at $a=0$ and iterating proves uniqueness on $[0,T_*]$.

We next establish nonnegativity.
Let $u_-:=\max\{-u,0\}$.
Testing the first equation by $u_-$, integrating by parts and using Young's inequality gives
\[
 \frac{1}{2}\frac{\dd}{\dd t}\norm{u_-}_{L^2}^2
 +\norm{\nabla_gu_-}_{L^2}^2
 \le\frac{1}{2}\norm{\nabla_gu_-}_{L^2}^2
    +\frac{1}{2}\norm{\nabla_gv}_{L^\infty}^2\norm{u_-}_{L^2}^2
 \quad\text{for all }0<t<T_{\max}.
\]
Since $\nabla_gv$ is bounded on $[0,T]$ by \eqref{eq:local-candidate-bounds}, Gronwall's inequality and $u_0\ge0$ imply that $u\ge0$.
Then \eqref{eq:local-w-v} and $w_0\ge0$ yield $w\ge0$.
Finally, \eqref{eq:kdef}, \eqref{eq:mass}, and \eqref{eq:local-w-v} yield \eqref{eq:local-mass-k}.

\medskip\noindent\textit{Step 4. Blow-up criterion.}
Steps 1--3 remain valid for nonnegative data in $C(\Sigma)\times C^\gamma(\Sigma)$, $\gamma\in(0,1)$, upon choosing $\alpha\in(0,\min\{\gamma,1-\frac{2}{p}\})$ in Step 2.
Since $w(a)\in C^\alpha(\Sigma)$ for every $a>0$ and the local existence time depends only on the $L^\infty$ norms of the initial data, the solution extends to a maximal interval $[0,T_{\max})$.
If $T_{\max}<\infty$ and
$$K:=\sup_{t\in[0,T_{\max})}\norm{u(t)}_{L^\infty}<\infty,$$ then
\eqref{eq:local-w-v} gives
\[
 \sup_{t\in[0,T_{\max})}\norm{w(t)}_{L^\infty}
 \le\norm{w_0}_{L^\infty}+K.
\]
The local construction gives an existence time $\tau>0$ independent of $a\in(0,T_{\max})$.
Restarting at $a>T_{\max}-\tau$ and using uniqueness extends the solution beyond $T_{\max}$, a contradiction.
Thus \eqref{eq:local-blow-up-alternative} holds.

\end{proof}

\section{Global existence}
\label{sec:global-existence}
In this section, we establish global existence using energy estimates and heat semigroup estimates.

\begin{lemma}\label{lem:finite-time-energy}
Let $(u,v,w)$ be the solution given by Proposition~\ref{prop:local-existence}.
There is a constant $C>0$, depending only on $(\Sigma,g)$ and the initial data, such that
\begin{equation*}
 Y(t)+1\le (Y(0)+1)e^{Ct}
 \quad \text{for all } t\in(0,T_{\max}),
\end{equation*}
where
\begin{equation*}
 Y(t):=\norm{u(t)}_{L^2}^2+\norm{w(t)}_{L^3}^3.
\end{equation*}
\end{lemma}

\begin{proof}
Multiplying the first and third equations in \eqref{eq:system} by $2u$ and $3w^2$, respectively, and integrating over $\Sigma$, we obtain, for all $t\in(0,T_{\max})$,
\begin{align*}
 \frac{\dd}{\dd t}\norm{u}_{L^2}^2+2\norm{\nabla_gu}_{L^2}^2
 &=-\overline w\norm{u}_{L^2}^2+\int_\Sigma u^2w\dd\mu_g,\\
 \frac{\dd}{\dd t}\norm{w}_{L^3}^3+3\norm{w}_{L^3}^3
 &=3\int_\Sigma uw^2\dd\mu_g.
\end{align*}
Adding these identities and using $\overline w\ge0$, we find
\begin{equation}\label{eq:Y-testing}
 Y'(t)+2\norm{\nabla_gu}_{L^2}^2+3\norm{w}_{L^3}^3
 \le\int_\Sigma u^2w\dd\mu_g+3\int_\Sigma uw^2\dd\mu_g
 \quad\text{for all }0<t<T_{\max}.
\end{equation}
For $\varepsilon>0$, Young's inequality implies that there is a constant $C_1(\varepsilon)>0$ such that
\begin{equation}\label{eq:mixed-terms-young}
 \int_\Sigma u^2w\dd\mu_g+3\int_\Sigma uw^2\dd\mu_g
 \le\varepsilon\norm{u}_{L^3}^3+C_1(\varepsilon)\norm{w}_{L^3}^3
 \quad\text{for all }0<t<T_{\max}.
\end{equation}
The Sobolev embedding $W^{1,1}(\Sigma)\hookrightarrow L^2(\Sigma)$ \cite[Theorem~2.6]{1999-Hebey}, H\"older's inequality, interpolation, and Young's inequality give
\begin{align*}
 \norm{u}_{L^3}^{\frac{3}{2}}
 &\le C_2\bigl(\norm{\nabla_g(u^{\frac{3}{2}})}_{L^1}
             +\norm{u^{\frac{3}{2}}}_{L^1}\bigr)\\
 &=\frac{3C_2}{2}\int_\Sigma u^{\frac{1}{2}}|\nabla_gu|\dd\mu_g
   +C_2\norm{u}_{L^{\frac{3}{2}}}^{\frac{3}{2}}\\
 &\le\frac{3C_2}{2}\norm{u^{\frac{1}{2}}}_{L^2}
                  \norm{\nabla_gu}_{L^2}
   +C_2\bigl(\norm{u}_{L^1}^{\frac{1}{2}}
            \norm{u}_{L^3}^{\frac{1}{2}}\bigr)^{\frac{3}{2}}\\
 &\le C_3\norm{u}_{L^1}^{\frac{1}{2}}\norm{\nabla_gu}_{L^2}
   +C_3\norm{u}_{L^1}^{\frac{3}{4}}\norm{u}_{L^3}^{\frac{3}{4}}\\
 &\le C_3\norm{u}_{L^1}^{\frac{1}{2}}\norm{\nabla_gu}_{L^2}
   +\frac12\norm{u}_{L^3}^{\frac{3}{2}}
   +\frac{C_3^2}{2}\norm{u}_{L^1}^{\frac{3}{2}}.
\end{align*}
Using mass conservation \eqref{eq:mass}, we obtain
\begin{equation}\label{eq:GN-L3-L1}
 \norm{u}_{L^3}^3
 \le C_4M\norm{\nabla_gu}_{L^2}^2
   +C_4M^3
 \quad\text{for all }0<t<T_{\max},
\end{equation}
where $C_2,C_3,C_4>0$ depend only on $(\Sigma,g)$.
Choose $\varepsilon>0$ with $\varepsilon C_4M\le1$.
Combining this estimate with \eqref{eq:Y-testing} and \eqref{eq:mixed-terms-young}, we have $Y'(t)\le C_5(Y(t)+1)$, where $C_5:=C_1(\varepsilon)+\varepsilon C_4M^3$.
Hence, using Gronwall's inequality and continuity at $t=0$, we complete the proof.
\end{proof}

We now prove Proposition~\ref{prop:global-existence}.

\begin{proof}[Proof of Proposition~\ref{prop:global-existence}]
Fix $T\in(0,T_{\max})$.
Lemma~\ref{lem:finite-time-energy} yields $C_1(T)>0$ such that
\begin{equation}\label{eq:global-L2-L3-bound}
 \sup_{t\in[0,T]}
 \bigl(\norm{u(t)}_{L^2}+\norm{w(t)}_{L^3}\bigr)\le C_1(T).
\end{equation}
By \eqref{eq:normalization}, elliptic regularity \cite[Theorem~A.3]{2019-Lee} and Sobolev embedding \cite[Theorem~2.7]{1999-Hebey} yield $C_2(T)>0$ such that
\begin{equation*}
 \sup_{t\in[0,T]}\norm{\nabla_gv(t)}_{L^\infty}
 \le C_2(T).
\end{equation*}
For $t\in(0,T]$, put $a:=\max\{0,t-1\}$.
Duhamel's formula gives
\begin{equation}\label{eq:global-duhamel}
 u(t)=e^{(t-a)\Delta_g}u(a)
 -\int_a^t e^{(t-s)\Delta_g}\operatorname{div}_g
       (u(s)\nabla_gv(s))\dd s.
\end{equation}
For $t\leqslant 1$, contractivity of the heat semigroup \cite[Exercise~7.33]{2009-Grigoryan} gives
$$\norm{e^{t\Delta_g}u_0}_{L^\infty}\le\norm{u_0}_{L^\infty}
 \quad\text{for all }0<t\le\min\{1,T\}.$$
For $t>1$, we have $a=t-1$, so Lemma~\ref{lem:semigroup-smoothing} and \eqref{eq:global-L2-L3-bound} yield $C_3>0$ such that
\[
 \norm{e^{(t-a)\Delta_g}u(a)}_{L^\infty}
 \le C_3\norm{u(a)}_{L^2}\le C_3C_1(T)
 \quad\text{for all }1<t\le T.
\]
Set $C_4(T):=\max\{\norm{u_0}_{L^\infty},C_3C_1(T)\}$.
By Lemma~\ref{lem:semigroup-smoothing}, \eqref{eq:global-L2-L3-bound}, \eqref{eq:global-duhamel}, and interpolation, we obtain
\begin{align*}
 \norm{u(t)}_{L^\infty}
 &\le \norm{e^{(t-a)\Delta_g}u(a)}_{L^\infty}
       +C_5\int_a^t(t-s)^{-\frac34}
       \norm{u(s)\nabla_gv(s)}_{L^4}\dd s\\
 &\le C_4(T)+C_5C_2(T)\int_a^t(t-s)^{-\frac34}
       \norm{u(s)}_{L^\infty}^{\frac12}
       \norm{u(s)}_{L^2}^{\frac12}\dd s\\
 &\le C_4(T)+C_5C_2(T)C_1(T)^{\frac12}
       \left(\sup_{r\in[0,T]}\norm{u(r)}_{L^\infty}\right)^{\frac12}
       \int_a^t(t-s)^{-\frac34}\dd s\\
 &\le C_4(T)+C_6(T)
       \left(\sup_{r\in[0,T]}\norm{u(r)}_{L^\infty}\right)^{\frac12}
 \quad\text{for all }0<t\le T,
\end{align*}
where $C_5,C_6(T)>0$ and we used $t-a\le1$.
Taking the supremum and using Young's inequality gives
\[
 \sup_{t\in[0,T]}\norm{u(t)}_{L^\infty}\le 2C_4(T)+C_6(T)^2.
\]
Together with \eqref{eq:local-blow-up-alternative}, this implies $T_{\max}=\infty$.
\end{proof}

\section{Uniform boundedness for subcritical mass}
\label{sec:subcritical}
In this section, we establish uniform boundedness for solutions with subcritical mass.
Let $A:=\mu_g(\Sigma)$.
We introduce an energy functional
\begin{equation}\label{eq:energy}
 \mathcal E_g(u,v)
 :=\int_\Sigma\bigl(L(u)-uv\bigr)\dd\mu_g
 +\frac{1}{2}\int_\Sigma|\nabla_gv|^2\dd\mu_g,
\end{equation}
where
\begin{equation}\label{eq:entropy-density}
 L(r):=r\log r-r+1\ge0\quad\text{for all }r\ge0,
\end{equation}
and write $\mathcal E_g(t):=\mathcal E_g(u(t),v(t))$.

\begin{lemma}\label{lem:energy-monotonicity}
Let $(u,v,w)$ be a global classical solution of \eqref{eq:system} with nonnegative initial data.
Then
\[
 \mathcal E_g(t)\le\mathcal E_g(0)\quad\text{for all }t\ge0.
\]
\end{lemma}

\begin{proof}
If $M=0$, then $u\equiv0$ and $\partial_t w=-w$.
The second equation in \eqref{eq:system} and \eqref{eq:normalization} imply $\partial_t v=-v$, so $\mathcal E_g'(t)=-\norm{\nabla_gv(t)}_{L^2}^2\le0$ for $t>0$.
Thus $\mathcal E_g(t)\le\mathcal E_g(0)$ by continuity at $t=0$.
Now assume $M>0$.
For every $\tau>0$, the parabolic strong maximum principle \cite[Proposition 2.3]{2006-Mueller} yields $u>0$ on $\Sigma\times[\tau,\infty)$.
Thus, for $t>\tau$, integration by parts on $\Sigma$ gives
\begin{align*}
 \frac{\dd}{\dd t}\mathcal E_g(t)
 &=\int_{\Sigma}(\log u-v)\partial_t u\dd\mu_g
   -\int_{\Sigma}u\partial_t v\dd\mu_g
   +\int_{\Sigma}\langle\nabla_gv,\nabla_g\partial_t v\rangle_g\dd\mu_g\\
 &=-\int_{\Sigma}u\left|\nabla_g(\log u-v)\right|^2\dd\mu_g
   -\int_{\Sigma}(u+\Delta_gv)\partial_t v\dd\mu_g.
\end{align*}
Since \eqref{eq:normalization} implies $\int_{\Sigma}\partial_t v\dd\mu_g=0$, using $\Delta_gv=\overline w-w$ and $\partial_t w=u-w$, we have
\begin{equation*}
\int_{\Sigma}(u+\Delta_gv)\partial_t v\dd\mu_g
=\int_{\Sigma}(\partial_t w+\overline w)\partial_t v\dd\mu_g
= \int_{\Sigma}(\partial_t w)(\partial_t v)\dd\mu_g.
\end{equation*}
Differentiating the second equation in \eqref{eq:system} in time yields
\[
 -\Delta_g\partial_t v=\partial_t w-\partial_t\overline w.
\]
Consequently,
\begin{align*}
 \int_{\Sigma}(\partial_t w)(\partial_t v)\dd\mu_g
 &=\int_{\Sigma}\bigl(-\Delta_g\partial_t v+\partial_t\overline w\bigr)\partial_t v\dd\mu_g
 =\int_{\Sigma}|\nabla_g\partial_t v|^2\dd\mu_g,
\end{align*}
and therefore
\begin{equation}\label{eq:dissipation}
 \frac{\dd}{\dd t}\mathcal E_g(t)
 =-\int_{\Sigma}u\left|\nabla_g(\log u-v)\right|^2\dd\mu_g
  -\int_{\Sigma}|\nabla_g\partial_t v|^2\dd\mu_g
 \le0
 \quad\text{for a.e. }t>0,
\end{equation}
which completes the proof.

\end{proof}

We now use this functional to derive the estimates needed when $M<8\pi$.

\begin{lemma}\label{lem:energy-entropy}
Let $u_0\in C(\Sigma)$ and $w_0\in C^1(\Sigma)$ be nonnegative, and let $(u,v,w)$ be the corresponding global classical solution of \eqref{eq:system}--\eqref{eq:normalization}.
If $$M=\int_\Sigma u_0\dd\mu_g<8\pi,$$ then there exists $C>0$ such that
\begin{equation*}
 \sup_{t\ge0}\int_{\Sigma}u(x,t)|\log u(x,t)|\dd\mu_g\le C.
\end{equation*}
\end{lemma}

\begin{proof}
If $M=0$, then $u\equiv0$ and the assertion is immediate.
Hence assume $M>0$.
Using the elementary inequality
\[
 a\log\frac{a}{b}-a+b\ge0
 \quad\text{for all }a\ge0,\ b>0,
\]
we take
\[
 a=u,
 \qquad
 b=\frac{Me^v}{\int_\Sigma e^v\dd\mu_g}.
\]
Integrating over $\Sigma$ and using \eqref{eq:mass}, we obtain
\[
 \begin{aligned}
 0&\le\int_\Sigma u\log\left(
 \frac{u\int_\Sigma e^v\dd\mu_g}{Me^v}\right)\dd\mu_g\\
 &=\int_\Sigma(u\log u-uv)\dd\mu_g
 -M\log M+M\log\left(\int_\Sigma e^v\dd\mu_g\right)
 \qquad\text{for all }t\ge0.
 \end{aligned}
\]
Hence
\begin{equation}\label{eq:entropy-variational}
 \int_\Sigma u\log u\dd\mu_g-\int_\Sigma uv\dd\mu_g
 \ge M\log M-M\log\left(\int_\Sigma e^v\dd\mu_g\right)
 \quad\text{for all }t\ge0.
\end{equation}
Applying the Moser-Trudinger inequality (see \cite{2007-TAMS-Yang}) yields
\begin{equation}\label{eq:exponential-v-bound}
 \int_\Sigma e^{|v|}\dd\mu_g
 \le C_1\exp\left(\frac{1}{16\pi}
 \|\nabla_gv\|_{L^2}^2\right)
 \quad\text{for all }t\ge0.
\end{equation}
Combining this with \eqref{eq:energy} and \eqref{eq:entropy-variational} yields
\begin{align*}
\mathcal E_g(t)
&=\int_{\Sigma}u\log u\dd\mu_g-\int_{\Sigma}uv\dd\mu_g-\int_{\Sigma}u\dd\mu_g+A+\frac{1}{2}\int_{\Sigma}\left|\nabla_gv\right|^2\dd\mu_g
\\
&\geq M\log M-M\log\left(\int_{\Sigma}e^v\dd\mu_g\right)
-M+A+\frac{1}{2}\int_{\Sigma}\left|\nabla_gv\right|^2\dd\mu_g
\\
 &\ge M\log\frac{M}{{C_1}}-M+A
 +\frac{8\pi-M}{16\pi}
   \int_{\Sigma}|\nabla_gv|^2\dd\mu_g
 \quad\text{for all }t\ge0.
\end{align*}
Together with $M<8\pi$ and Lemma~\ref{lem:energy-monotonicity}, this implies that there exists ${C_2}>0$ such that
\begin{equation}\label{eq:vH1}
 \sup_{t\ge0}\norm{\nabla_gv(t)}_{L^2}\le {C_2}.
\end{equation}
The convex conjugate of $L$ is $L^*(s)=e^s-1$.
Hence Fenchel's inequality \cite[p.~94]{2004-BoydVandenberghe} gives
\begin{equation*}
 2uv\le L(u)+e^{2v}-1
 \quad\text{for all }x\in\Sigma,\ t\ge0.
\end{equation*}
Using Lemma~\ref{lem:energy-monotonicity} together with this pointwise bound and \eqref{eq:energy}, we obtain
\begin{align*}
 \int_{\Sigma}L(u)\dd\mu_g
 &=\mathcal E_g(t)+\int_{\Sigma}uv\dd\mu_g
 -\frac{1}{2}\int_{\Sigma}|\nabla_gv|^2\dd\mu_g\\
 &\le\mathcal E_g(0)+\int_{\Sigma}uv\dd\mu_g\\
 &\le \mathcal E_g(0)
 +\frac{1}{2}\int_{\Sigma}L(u)\dd\mu_g
 +\frac{1}{2}\int_{\Sigma}(e^{2v}-1)\dd\mu_g
 \quad\text{for all }t\ge0,
\end{align*}
whence
\begin{equation*}
 \int_{\Sigma}L(u)\dd\mu_g
 \le 2\mathcal E_g(0)+\int_{\Sigma}e^{2v}\dd\mu_g
 \quad\text{for all }t\ge0.
\end{equation*}
Then, applying \eqref{eq:exponential-v-bound} to $2v$ and using \eqref{eq:vH1}, we find
\[
 \sup_{t\ge0}\int_{\Sigma}L(u(x,t))\dd\mu_g<\infty.
\]
Since $s\log s\le L(s)+s$ for $s\ge1$ and $s|\log s|\le e^{-1}$ for $0<s<1$, \eqref{eq:mass} gives
\[
\begin{aligned}
\int_\Sigma u|\log u|\,\mathrm d\mu_g
&\le \int_{\{u\ge1\}}\bigl(L(u)+u\bigr)\,\mathrm d\mu_g
+\frac{\mu_g(\{u<1\})}{e}\\
&\le \int_\Sigma L(u)\,\mathrm d\mu_g+M+\frac Ae.
\end{aligned}
\]
This completes the proof.
\end{proof}

We next establish an interpolation inequality that will be used to derive an $L^2$ bound for $u$.
\begin{lemma}\label{lem:entropy-interpolation}
Let $B>0$.
For every $\varepsilon>0$, there exists $C_\varepsilon>0$ such that
\begin{equation*}
 \norm{f}_{L^3}^3
 \le \varepsilon\norm{\nabla_gf}_{L^2}^2+C_\varepsilon
\end{equation*}
for all nonnegative $f\in W^{1,2}(\Sigma)$ satisfying
\begin{equation*}
 \int_\Sigma f|\log f|\dd\mu_g\le B.
\end{equation*}
\end{lemma}

\begin{proof}
Fix $K>e$ and let $f_K:=(f-K)_+$.
Then
\begin{equation}\label{eq:tail-L1}
 \norm{f_K}_{L^1}
 \le \int_{\{f>K\}}f\dd\mu_g
\le\frac{1}{\log K}\int_\Sigma f|\log f|\dd\mu_g
 \le\frac{B}{\log K}.
\end{equation}
By \eqref{eq:GN-L3-L1}, there exists $C_1>0$ such that
\[
 \norm{f_K}_{L^3}^3
 \le C_1\norm{\nabla_gf_K}_{L^2}^2\norm{f_K}_{L^1}
 +C_1\norm{f_K}_{L^1}^3.
\]
Using $f\le K+f_K$, \eqref{eq:tail-L1}, and $|\nabla_gf_K|\le|\nabla_gf|$ a.e., we obtain
\begin{align*}
 \norm{f}_{L^3}^3 
 &=\int_{\Sigma}f^3\dd\mu_g \\
 &\leq 4\int_{\Sigma}f_K^3\dd\mu_g + 4\int_{\Sigma}K^3\dd\mu_g\\
 &= 4\norm{f_K}_{L^3}^3+4 AK^3\\
 &\leq  \frac{4{C_1}B}{\log K}\left\|\nabla_g f\right\|_{L^2}^2+4{C_1}\left(\frac{B}{\log K}\right)^3+4 AK^3.
\end{align*}
Choosing $K$ sufficiently large so that $\frac{4C_1B}{\log K}\le\varepsilon$ completes the proof.
\end{proof}

\begin{lemma}\label{lem:L2L3}
Under the assumptions of Lemma~\ref{lem:energy-entropy}, there exists $C>0$ such that
\begin{equation*}
 \sup_{t\ge0}\left(
 \norm{u(t)}_{L^2}
 +\norm{w(t)}_{L^3}
 \right)\le C.
\end{equation*}
\end{lemma}

\begin{proof}
Let $Y$ be as defined in Lemma~\ref{lem:finite-time-energy}.
Since Young's inequality implies
\begin{equation*}
 \int_\Sigma u^2w\dd\mu_g+3\int_\Sigma uw^2\dd\mu_g
 \le\frac{5}{3}\norm{u}_{L^3}^3+\frac{7}{3}\norm{w}_{L^3}^3
 \quad\text{for all }t>0,
\end{equation*}
using \eqref{eq:Y-testing}, we have
\begin{equation}\label{eq:Y-pre}
 Y'(t)+2\norm{\nabla_gu}_{L^2}^2+\frac{2}{3}\norm{w}_{L^3}^3
 \le\frac{5}{3}\norm{u}_{L^3}^3
 \quad\text{for all }t>0.
\end{equation}
Combining Lemma~\ref{lem:energy-entropy} with Lemma~\ref{lem:entropy-interpolation} applied with $\varepsilon=\frac{3}{5}$, we obtain
\[
 \frac{5}{3}\norm{u}_{L^3}^3
 \le \norm{\nabla_gu}_{L^2}^2+{C_1}
 \quad\text{for all }t>0.
\]
Substituting this into \eqref{eq:Y-pre}, we obtain
\begin{equation}\label{eq:Y-diss}
 Y'(t)+\norm{\nabla_gu}_{L^2}^2+\frac{2}{3}\norm{w}_{L^3}^3\le {C_1}
 \quad\text{for all }t>0.
\end{equation}
By the Sobolev--Poincar\'e inequality \cite[Theorem~2.11]{1999-Hebey}, H\"older's inequality, and \eqref{eq:mass}, there exists $C_2>0$ such that
\begin{equation*}
 \norm{u}_{L^2}^2
 =\left\|u-\frac{M}{A}\right\|_{L^2}^2+\frac{M^2}{A}
 \le C_2\norm{\nabla_gu}_{L^2}^2+\frac{M^2}{A}
 \quad\text{for all }t>0.
\end{equation*}
Combining this with  \eqref{eq:Y-diss}, we obtain ${C_3},{C_4}>0$ such that
\[
 Y'(t)+C_3Y(t)\le {C_4}
 \quad\text{for all }t>0,
\]
and Gronwall's inequality yields the asserted bound.
\end{proof}

We are now ready to prove Theorem~\ref{thm:subcritical}.

\begin{proof}[Proof of Theorem~\ref{thm:subcritical}]
By Lemma~\ref{lem:L2L3}, there exists $C_1>0$ such that
\begin{equation*}
 \sup_{t\ge0}\norm{w(t)}_{L^3}\le C_1.
\end{equation*}
By \eqref{eq:normalization}, elliptic regularity \cite[Theorem~A.3]{2019-Lee} and the Sobolev embedding \cite[Theorem~2.7]{1999-Hebey}, there exist $C_2,C_3>0$ such that
\begin{equation}\label{eq:grad-v-infty}
\norm{\nabla_gv(t)}_{L^\infty}\le C_2
\norm{v(t)}_{W^{2,3}}\le
C_3 \quad \text{for all } t\ge0.
\end{equation}
For $t\ge1$, Duhamel's formula on $[t-1,t]$ reads
\begin{equation*}
 u(t)=e^{\Delta_g}u(t-1)
 -\int_{t-1}^{t}e^{(t-s)\Delta_g}
  \operatorname{div}_g\bigl(u(s)\nabla_gv(s)\bigr)\dd s.
\end{equation*}
Fix $2<p<\infty$.
Applying Lemmas~\ref{lem:semigroup-smoothing} and~\ref{lem:L2L3} together with \eqref{eq:grad-v-infty} gives ${C_4},{C_5},{C_6}>0$ such that
\begin{align*}
 \norm{u(t)}_{L^p}
 &\le {C_4}\norm{u(t-1)}_{L^2}
 +{C_4}\int_{t-1}^{t}(t-s)^{-1+\frac{1}{p}}
   \norm{u(s)\nabla_gv(s)}_{L^2}\dd s\\
 &\le {C_5}+{C_5}\int_0^1\tau^{-1+\frac{1}{p}}\dd\tau\le {C_6}
 \quad \text{for all } t\ge1.
\end{align*}
Using this and Lemma~\ref{lem:semigroup-smoothing} in the same argument, we obtain $C_7,C_8,C_9>0$ such that
\begin{align*}
 \norm{u(t)}_{L^\infty}
 &\le {C_7}\norm{u(t-1)}_{L^p}
 +{C_7}\int_{t-1}^{t}(t-s)^{-\frac{1}{2}-\frac{1}{p}}
   \norm{u(s)\nabla_gv(s)}_{L^p}\dd s\\
 &\le {C_8}+{C_8}\int_0^1\tau^{-\frac{1}{2}-\frac{1}{p}}\dd\tau\le {C_9}
  \quad \text{for all } t\ge2.
\end{align*}
The continuity of the solution on the compact interval $[0,2]$ then yields
\begin{equation}\label{eq:u-infty}
 \sup_{t\ge0}\norm{u(t)}_{L^\infty}\le {C_{10}}
\end{equation}
for some ${C_{10}}>0$.
Finally, the third equation in \eqref{eq:system} gives
\begin{equation}\label{eq:w-representation}
 w(t)=e^{-t}w_0+\int_0^t e^{-(t-s)}u(s)\dd s.
\end{equation}
Combining this with \eqref{eq:u-infty}, we obtain $C_{11}>0$ such that
\[
 \sup_{t\ge0}\norm{w(t)}_{L^\infty}\le C_{11}.
\]
Together with \eqref{eq:grad-v-infty}, this completes the proof.
\end{proof}

\section{Unbounded solutions for supercritical mass}
\label{sec:supercritical}

For $M>0$, let $\mathcal S_M$ denote the set of all smooth stationary solutions $(u,v,w)$ of \eqref{eq:system} satisfying $u,w\ge0$, $\int_\Sigma u\dd\mu_g=M$, and $\int_\Sigma v\dd\mu_g=0$.
This set is nonempty, since $\left(\frac{M}{A},0,\frac{M}{A}\right)\in\mathcal S_M$.
Recalling  \eqref{eq:mass} and $A=\mu_g(\Sigma)$, we see that
\begin{equation}\label{eq:stationary-triple}
 u=w=\frac{Me^{v}}{\int_\Sigma e^{v}\dd\mu_g}.
\end{equation}
and
\begin{equation}\label{eq:stationary-mean-field}
 -\Delta_gv
 =M\left(\frac{e^{v}}{\int_\Sigma e^{v}\dd\mu_g}-\frac{1}{A}\right),
 \quad \int_\Sigma v\dd\mu_g=0.
\end{equation}

We first establish an estimate for a weighted mean field equation.
\begin{lemma}\label{lem:weighted-energy}
Assume that $\Sigma$ is orientable.
Fix $\lambda>0$, $K\ge1$, and $\kappa>0$.
If $V,W\in W^{1,\infty}(\Sigma)$ satisfy
\begin{equation}\label{eq:coefficient-bounds}
V\ge\kappa,\qquad
\norm V_{W^{1,\infty}(\Sigma)}+\norm W_{W^{1,\infty}(\Sigma)}\le K,
\qquad \int_\Sigma W\dd\mu_g=1,
\end{equation}
and $\eta\in C^2(\Sigma)$ satisfies
\begin{equation}\label{eq:closed-equation}
-\Delta_g\eta=\lambda(Ve^\eta-W),
\qquad \int_\Sigma Ve^\eta\dd\mu_g=1,
\end{equation}
then there is a constant $C=C(\Sigma,g,\lambda,K,\kappa)$ such that
\begin{equation*}
\int_\Sigma\eta Ve^\eta\dd\mu_g+\overline\eta\ge-C.
\end{equation*}
\end{lemma}

\begin{proof}
Suppose that the assertion fails.
Then there are sequences $V_j,W_j,\eta_j$ satisfying the hypotheses such that
\begin{equation}\label{eq:bad-energy}
\int_\Sigma\eta_jV_je^{\eta_j}\dd\mu_g+\overline\eta_{j,g}
\longrightarrow-\infty.
\end{equation}
By \eqref{eq:closed-equation}, the absolute value of the left-hand side of \eqref{eq:bad-energy} is bounded by $2\norm{\eta_j}_{L^\infty(\Sigma)}$.
Hence $\norm{\eta_j}_{L^\infty(\Sigma)}\to\infty$.
Applying \cite[Theorem~0.2]{1999-CMP-Li} and passing to a subsequence, we obtain distinct points $p_1,\ldots,p_N$ and local maximum points $p_{j,\ell}\to p_\ell$, with $h_{j,\ell}:=\eta_j(p_{j,\ell})\to\infty$.
As in \cite[Theorem~0.2(c)]{1999-CMP-Li}, let $y$ be an isothermal coordinate centered at $p_{j,\ell}$.
For fixed sufficiently small $r>0$ and all large $j$, the disks $\mathcal U_{j,\ell}:=\left\{x\in\Sigma:d(x,p_{j,\ell})<r\right\}$ are disjoint and contain $p_\ell$.
The estimates in \cite[Theorem~0.2(c)--(e)]{1999-CMP-Li} give constants $C_1,C_2,C_3>0$, independent of $j$, such that
\begin{align}
\left|\eta_j(y)-h_{j,\ell}
+2\log\left(1+c_{j,\ell}e^{h_{j,\ell}}|y|^2\right)\right|
&\le C_1,\label{eq:bubble-profile}\\
|h_{j,\ell}+\overline\eta_{j,g}|&\le C_2,
\label{eq:local-height-mean}
\end{align}
where $c_{j,\ell}:=\lambda V_j(p_{j,\ell})/8$ is bounded above and away from zero, and
\begin{equation}\label{eq:outside-estimate}
|\eta_j-\overline\eta_{j,g}|\le C_3
\quad\text{on }\Sigma\setminus\bigcup_{\ell=1}^N\mathcal U_{j,\ell}.
\end{equation}
Let $H_j:=\max_\Sigma\eta_j$.
By \eqref{eq:local-height-mean} and \eqref{eq:outside-estimate}, outside the disks we have $\eta_j\le -h_{j,1}+C_2+C_3<h_{j,1}$ for large $j$.
Thus \cite[Theorem~0.2(b)]{1999-CMP-Li} gives $H_j=\max_{1\le \ell\le N} h_{j,\ell}$.
Hence \eqref{eq:local-height-mean} yields
\begin{equation}\label{eq:height-mean}
|H_j+\overline\eta_{j,g}|\le C_2,
\qquad 0\le H_j-h_{j,\ell}\le2C_2.
\end{equation}
Using \eqref{eq:coefficient-bounds}, \eqref{eq:bubble-profile} and \eqref{eq:height-mean}, and then changing variables to $z=c_{j,\ell}^{1/2}e^{h_{j,\ell}/2}y$, we find $C_4,C_5>0$ such that
\begin{equation}
\int_{\mathcal U_{j,\ell}}(H_j-\eta_j)V_je^{\eta_j}\dd\mu_g
\le C_4\int_{\mathbb R^2}
\frac{1+\log(1+|z|^2)}{(1+|z|^2)^2}\dd z \le C_5.
\label{eq:local-deficit}
\end{equation}
By \eqref{eq:outside-estimate} and \eqref{eq:height-mean}, we obtain $C_6>0$ such that
\begin{equation}\label{eq:outside-deficit}
\int_{\Sigma\setminus\bigcup_\ell\mathcal U_{j,\ell}}
(H_j-\eta_j)V_je^{\eta_j}\dd\mu_g
\le KA(2H_j+C_2+C_3)e^{C_2+C_3-H_j}\le C_6.
\end{equation}
Using \eqref{eq:height-mean}--\eqref{eq:outside-deficit}, we obtain
\begin{align*}
\int_\Sigma\eta_jV_je^{\eta_j}\dd\mu_g+\overline\eta_{j,g}
&=H_j+\overline\eta_{j,g}
-\int_\Sigma(H_j-\eta_j)V_je^{\eta_j}\dd\mu_g\\
&\ge-C_2-NC_5-C_6,
\end{align*}
contradicting \eqref{eq:bad-energy}.
This completes the proof.
\end{proof}

We now apply this estimate to the stationary solutions on $\Sigma$.

\begin{proposition}\label{prop:stationary-energy-lower}
Let $M>0$.
Then the energy functional $\mathcal E_g$ defined in \eqref{eq:energy} satisfies
\begin{equation}\label{eq:stationary-energy-lower-bound}
 \beta_M:=\inf_{(u,v,w)\in\mathcal S_M}
 \mathcal E_g(u,v)> -\infty.
\end{equation}
\end{proposition}

\begin{proof}
For $(u,v,w)\in\mathcal S_M$, set
\[
 Z:=\int_\Sigma e^{v}\dd\mu_g,
 \quad \xi:=v-\log Z,
 \quad \overline\xi:=\frac1A\int_\Sigma\xi\dd\mu_g.
\]
Since $v$ has zero mean, $v=\xi-\overline\xi$.
By \eqref{eq:stationary-mean-field},
\begin{equation}\label{eq:stationary-normalized}
 -\Delta_g\xi=M\left(e^{\xi}-\frac1A\right),
 \quad \int_\Sigma e^{\xi}\dd\mu_g=1.
\end{equation}
Moreover, \eqref{eq:stationary-triple} implies $u=Me^{\xi}$ and $\log u-v=\log M+\overline\xi$.
Substituting these identities into \eqref{eq:energy} and using \eqref{eq:stationary-normalized}, we obtain
\begin{align*}
 \mathcal E_g(u,v)
 &=M\log M-M+A+M\overline\xi
   +\frac12\norm{\nabla_g\xi}_{L^2}^2\\
 &=M\log M-M+A
   +\frac M2\left(\int_\Sigma\xi e^{\xi}\dd\mu_g
   +\overline\xi\right).
\end{align*}
If $\Sigma$ is orientable, Lemma~\ref{lem:weighted-energy} with $\lambda=M$, $V=1$, and $W=1/A$ yields \eqref{eq:stationary-energy-lower-bound}.
If $\Sigma$ is nonorientable, let $\pi:(\widetilde\Sigma,\widetilde g)\to(\Sigma,g)$ be its connected orientable double cover with metric \(\widetilde g=\pi^*g\) \cite[Theorem~1.58 and Definition~1.59]{1982-Aubin}.
The lift $(u\circ\pi,v\circ\pi,w\circ\pi)$ is a stationary triple on $\widetilde\Sigma$ of mass $2M$.
Define the energy $\widetilde{\mathcal E}_g$ on $\widetilde\Sigma$ by
\[
\widetilde{\mathcal E}_g(u\circ\pi,v\circ\pi):=2\mathcal E_g(u,v).
\]
Applying the orientable case to $\widetilde{\mathcal E}_g$ implies \eqref{eq:stationary-energy-lower-bound}.
\end{proof}

We next show that bounded solutions converge subsequentially to stationary solutions.

\begin{lemma}\label{lem:bounded-stationary-limit}
Let $u_0\in C(\Sigma)$ and $w_0\in C^1(\Sigma)$ be nonnegative, with $M>0$.
Suppose that the corresponding global classical solution of \eqref{eq:system} satisfies
\begin{equation}\label{eq:bounded-trajectory-assumption}
 \sup_{t\ge0}\norm{u(t)}_{L^\infty}\le K<\infty.
\end{equation}
Then there exist a sequence $(t_j)_{j\ge1}$ with $t_j\to\infty$ as $j\to\infty$ and a stationary solution $(u_*,v_*,w_*)\in\mathcal S_M$ such that
\begin{equation*}
 u(t_j)\to u_*\ \text{in }L^2(\Sigma),
 \quad v(t_j)\to v_*\ \text{in }C^1(\Sigma)
 \quad\text{as }j\to\infty,
\end{equation*}
and
\begin{equation}\label{eq:bounded-limit-energy}
 \mathcal E_g(u_*,v_*)=\lim_{j\to\infty}\mathcal E_g(t_j)
 \le\mathcal E_g(0),
\end{equation}
where $\mathcal E_g$ is defined in \eqref{eq:energy}.
\end{lemma}

\begin{proof}
The representation \eqref{eq:w-representation} and \eqref{eq:bounded-trajectory-assumption} give, for all $t>0$,
\[
 \norm{w(t)}_{L^\infty}
 \le e^{-t}\norm{w_0}_{L^\infty}+K(1-e^{-t})
 \le\norm{w_0}_{L^\infty}+K=:C_1
 \quad\text{for all }t\ge0.
\]
For $2<p<\infty$, elliptic regularity \cite[Theorem~A.3]{2019-Lee} and Sobolev embedding \cite[Theorem~2.7]{1999-Hebey} give $C_2,C_3>0$ such that
\begin{equation}
\norm{v(t)}_{L^\infty}
 +\norm{\nabla_gv(t)}_{L^\infty}
 \le C_2\norm{v(t)}_{W^{2,p}}\le C_3
 \quad \text{for all } t>0.
 \label{eq:bounded-trajectory-elliptic}
\end{equation}
The nonnegativity of $L$ in \eqref{eq:entropy-density} implies
\begin{equation}\label{eq:bounded-trajectory-energy-lower}
 \mathcal E_g(t)\ge-\int_\Sigma uv\dd\mu_g\ge-MC_3
 \quad\text{for all }t>0.
\end{equation}
Define
\[
 \mathcal D(t):=\int_\Sigma
 u|\nabla_g(\log u-v)|^2\dd\mu_g
 +\norm{\nabla_g\partial_t v}_{L^2}^2
 \quad \text{for all } t>0.
\]
Integration of \eqref{eq:dissipation}, together with \eqref{eq:bounded-trajectory-energy-lower}, gives
\begin{equation*}
 \int_1^\infty\mathcal D(t)\dd t
 \le\mathcal E_g(1)+MC_3<\infty.
\end{equation*}
Then, the nonnegativity of $\mathcal D$ implies
\[
 \int_j^{j+1}\mathcal D(t)\dd t
 \le  \int_j^{\infty}\mathcal D(t)\dd t
 \to0
 \quad\text{as }j\to\infty.
\]
Since $\mathcal D$ is continuous, we can choose a sequence $(t_j)_{j\ge1}$ with $t_j\in[j,j+1]$ for each $j\ge1$ such that
\[
 0\le\mathcal D(t_j)\le\int_j^{j+1}\mathcal D(t)\dd t
 \quad\text{for all }j\ge1,
\]
which yields
\begin{equation}\label{eq:vanishing-dissipation-sequence}
 \mathcal D(t_j)\to0  \quad\text{as }j\to\infty.
\end{equation}
The strong maximum principle \cite[Proposition 2.3]{2006-Mueller} gives $u(t_j)>0$.
Since $u(t_j)\le K$, we have
\begin{equation}\label{eq:bounded-u-dissipation}
 \norm{\nabla_gu(t_j)-u(t_j)\nabla_gv(t_j)}_{L^2}^2
 =\int_\Sigma u(t_j)^2|\nabla_g(\log u(t_j)-v(t_j))|^2\dd\mu_g
 \le K\mathcal D(t_j)
 \quad\text{for all }j\ge1.
\end{equation}
Combining this with \eqref{eq:bounded-trajectory-elliptic}, we obtain
\begin{align*}
 \norm{\nabla_gu(t_j)}_{L^2}
 &\le\norm{\nabla_gu(t_j)-u(t_j)\nabla_gv(t_j)}_{L^2}
 +\norm{u(t_j)\nabla_gv(t_j)}_{L^2}\\
 &\le K^{\frac12}\mathcal D(t_j)^{\frac12}
 +(KM)^{\frac12}C_3
 \quad\text{for all }j\ge1.
\end{align*}
The compact embeddings in \cite[Theorem~2.9]{1999-Hebey} therefore give, along a subsequence,
\begin{equation}
 u(t_j)\rightharpoonup u_*\ \text{in }W^{1,2}(\Sigma),\quad
 u(t_j)\to u_*\ \text{in }L^2(\Sigma)\text{ and a.e.}\quad
 \quad\text{as }j\to\infty,
 \label{eq:u-sequence-convergence}
\end{equation}
and
\begin{equation}\label{eq:v-sequence-convergence}
 v(t_j)\to v_*\ \text{in }C^1(\Sigma)
 \quad\text{as }j\to\infty.
\end{equation}
The limits satisfy $0\le u_*\le K$ and
\begin{equation}\label{eq:limit-mass-normalization}
 \int_\Sigma u_*\dd\mu_g=M,
 \qquad \int_\Sigma v_*\dd\mu_g=0.
\end{equation}
Let $w_*=u_*$.
It remains to verify that $(u_*,v_*,w_*)\in\mathcal S_M$.
By \eqref{eq:vanishing-dissipation-sequence}--\eqref{eq:v-sequence-convergence}, weak lower semicontinuity of the $L^2$ norm yields
\[
 \nabla_gu_*=u_*\nabla_gv_*.
\]
Hence $\nabla_g(u_*e^{-v_*})=0$.
Since $\Sigma$ is connected, the Poincar\'e inequality \cite[Theorem~2.11]{1999-Hebey} gives $u_*e^{-v_*}=C_4$ almost everywhere for some $C_4\ge0$.
By \eqref{eq:limit-mass-normalization},
\begin{equation}\label{eq:limit-boltzmann-relation}
 u_* =C_4e^{v_*}
 =\frac{Me^{v_*}}{\int_\Sigma e^{v_*}\dd\mu_g}.
\end{equation}
Combining the second and third equations in \eqref{eq:system} with \eqref{eq:mass}, we obtain
\begin{equation}\label{eq:potential-relaxation-identity}
 -\Delta_g(\partial_t v+v)=u-\frac{M}{A}.
\end{equation}
Since $\int_\Sigma \partial_t v\dd\mu_g=0$, the Poincar\'e inequality \cite[Theorem~2.11]{1999-Hebey} and \eqref{eq:vanishing-dissipation-sequence} give $C_5>0$ such that
\begin{equation}\label{eq:limit-vt-vanishing}
 \norm{\partial_t v(t_j)}_{W^{1,2}}
 \le C_5\norm{\nabla_g\partial_t v(t_j)}_{L^2}\to0
 \quad\text{as }j\to\infty.
\end{equation}
Testing \eqref{eq:potential-relaxation-identity} by $\zeta\in C^\infty(\Sigma)$ yields
\[
 \int_\Sigma\langle\nabla_g(v(t_j)+\partial_t v(t_j)),\nabla_g\zeta\rangle_g
 \dd\mu_g
 =\int_\Sigma\left(u(t_j)-\frac{M}{A}\right)\zeta\dd\mu_g.
\]
Passing to the limit as $j\to\infty$ using \eqref{eq:u-sequence-convergence}, \eqref{eq:v-sequence-convergence} and \eqref{eq:limit-vt-vanishing} gives
\[
 -\Delta_gv_*=u_*-\frac{M}{A}
 \quad\text{in the sense of distributions}.
\]
Together with \eqref{eq:limit-boltzmann-relation}, this yields \eqref{eq:stationary-mean-field}.
Since $v_*\in C^1(\Sigma)$ yields $u_*\in C^1(\Sigma)$ by \eqref{eq:limit-boltzmann-relation}, elliptic regularity \cite[Theorem~A.6]{2019-Lee} gives $v_*\in C^{2+\alpha}(\Sigma)$ for every \(\alpha\in(0,1)\) and, by iteration, $v_*,u_*\in C^\infty(\Sigma)$.
Thus $(u_*,v_*,w_*)\in\mathcal S_M$.

Finally, consider the energy convergence in \eqref{eq:bounded-limit-energy}.
Since $L$ is continuous on $[0,K]$, \eqref{eq:u-sequence-convergence} and dominated convergence give $\int_\Sigma L(u(t_j))\dd\mu_g\to\int_\Sigma L(u_*)\dd\mu_g$ as $j\to\infty$.
By \eqref{eq:bounded-trajectory-elliptic}, \eqref{eq:u-sequence-convergence} and \eqref{eq:v-sequence-convergence}, the mixed term satisfies
\[
 \left|\int_\Sigma(u(t_j)v(t_j)-u_*v_*)\dd\mu_g\right|
 \le C_3\norm{u(t_j)-u_*}_{L^1}
 +M\norm{v(t_j)-v_*}_{L^\infty}\to0
 \quad\text{as }j\to\infty.
\]
By \eqref{eq:v-sequence-convergence}, we also have $\norm{\nabla_gv(t_j)}_{L^2}^2\to\norm{\nabla_gv_*}_{L^2}^2$ as $j\to\infty$.
Thus $\mathcal E_g(u(t_j),v(t_j))\to\mathcal E_g(u_*,v_*)$ as $j\to\infty$, and Lemma~\ref{lem:energy-monotonicity} yields \eqref{eq:bounded-limit-energy}.
\end{proof}

We then construct initial data with arbitrarily negative energy for $M>8\pi$.

\begin{lemma}\label{lem:supercritical-low-energy}
For every $M>8\pi$ and $R\in\mathbb R$, there exist strictly positive $u_0,w_0\in C^\infty(\Sigma)$ satisfying $\int_\Sigma u_0\dd\mu_g=M$ and
\[
 \mathcal E_g(u_0,v_0)<R,
\]
where $\mathcal E_g$ is the energy functional defined in \eqref{eq:energy}, and $v_0$ is determined by
\[
 -\Delta_gv_0
 =w_0-\frac{1}{A}\int_\Sigma w_0\dd\mu_g,
 \quad \int_\Sigma v_0\dd\mu_g=0.
\]
\end{lemma}

\begin{proof}
Fix $x_0\in\Sigma$ and choose $r_0>0$ with $2r_0$ smaller than the injectivity radius at $x_0$.
Set $r(x):=d_g(x,x_0)$ and choose a smooth cutoff $0\le\chi\le1$ such that $\chi=1$ on $B_g(x_0,r_0)$ and $\chi=0$ outside $B_g(x_0,2r_0)$.
Define
\[
 \psi(x):=\chi(x)r(x)^2+(1-\chi(x))r_0^2.
\]
Then $\psi\in C^\infty(\Sigma)$, $\psi=r^2$ on $B_g(x_0,r_0)$, and $\psi\ge r_0^2$ outside $B_g(x_0,r_0)$.
For $0<\varepsilon<\min\{1,r_0\}$, put
\begin{equation*}
 b_\varepsilon:=-2\log(\varepsilon^2+\psi),\quad
 \overline b_\varepsilon:=\frac{1}{A}\int_\Sigma b_\varepsilon\dd\mu_g,
 \quad \phi_\varepsilon:=b_\varepsilon-\overline b_\varepsilon,
\end{equation*}
and
\begin{equation*}
 u_{0,\varepsilon}
 :=\frac{Me^{\phi_\varepsilon}}
         {\int_\Sigma e^{\phi_\varepsilon}\dd\mu_g},
 \quad
 w_{0,\varepsilon}:=c_\varepsilon-\Delta_g\phi_\varepsilon,
 \quad c_\varepsilon:=1+\max\{0,\max_\Sigma\Delta_g\phi_\varepsilon\}.
\end{equation*}
Notice that $u_{0,\varepsilon}>0$, $w_{0,\varepsilon}\ge1$ and $\int_\Sigma u_{0,\varepsilon}\dd\mu_g=M$.
Since $\int_\Sigma\Delta_g\phi_\varepsilon\dd\mu_g=0$ by the divergence theorem, we have
\[
 -\Delta_g\phi_\varepsilon
 =w_{0,\varepsilon}-\frac{1}{A}\int_\Sigma w_{0,\varepsilon}\dd\mu_g,
 \quad \int_\Sigma\phi_\varepsilon\dd\mu_g=0.
\]
We therefore set $v_{0,\varepsilon}:=\phi_\varepsilon$.
The definition of $u_{0,\varepsilon}$ gives
\[
 \log u_{0,\varepsilon}
 =\log M+\phi_\varepsilon
   -\log\left(\int_\Sigma e^{\phi_\varepsilon}\dd\mu_g\right).
\]
Substituting this into \eqref{eq:energy} and using $\int_\Sigma u_{0,\varepsilon}\dd\mu_g=M$, we obtain
\begin{equation}\label{eq:bubble-initial-energy}
 \mathcal E_g(u_{0,\varepsilon},v_{0,\varepsilon})
 =\frac{1}{2}\norm{\nabla_g\phi_\varepsilon}_{L^2}^2
   -M\log\left(\int_\Sigma e^{\phi_\varepsilon}\dd\mu_g\right)
   +M\log M-M+A.
\end{equation}
In geodesic polar coordinates on $B_g(x_0,r_0)$, the volume element is
\begin{equation}\label{eq:polar-volume-expansion}
 \dd\mu_g=r(1+O(r^2))\dd r\mathrm{d}\theta
 \quad \text{for all } 0\le\theta<2\pi.
\end{equation}
Since $\varepsilon<\min\{1,r_0\}$, we have $r^2\le\varepsilon^2+r^2\le2r_0^2$ on $B_g(x_0,r_0)$ and $r_0^2\le\varepsilon^2+\psi\le1+\max_\Sigma\psi$ outside this ball.
Thus, by \eqref{eq:polar-volume-expansion}, there exist constants $C_1,C_2>0$, independent of $\varepsilon$, such that
\begin{align}
 |\overline b_\varepsilon|
 &\le\frac{1}{A}\int_\Sigma|b_\varepsilon|\dd\mu_g\notag\\
 &=\frac{2}{A}\int_{B_g(x_0,r_0)}
       |\log(\varepsilon^2+r^2)|\dd\mu_g
   +\frac{2}{A}\int_{\Sigma\setminus B_g(x_0,r_0)}
       |\log(\varepsilon^2+\psi)|\dd\mu_g\notag\\
 &\le\frac{1}{A}\int_{B_g(x_0,r_0)}
       (C_1+4|\log r|)\dd\mu_g
   +\frac{C_1}{A}\mu_g(\Sigma\setminus B_g(x_0,r_0))\notag\\
 &\le C_1\left(1+\int_0^{r_0}r|\log r|\dd r\right)\le C_2
 \quad\text{for all }0<\varepsilon<\min\{1,r_0\}.
 \label{eq:bubble-mean-bound}
\end{align}
By \eqref{eq:polar-volume-expansion} and $\psi\ge r_0^2$ outside $B_g(x_0,r_0)$,
\begin{align*}
 \int_\Sigma e^{b_\varepsilon}\dd\mu_g
 &=\int_{B_g(x_0,r_0)}\frac{\dd\mu_g}{(\varepsilon^2+r^2)^2}
   +\int_{\Sigma\setminus B_g(x_0,r_0)}
     \frac{\dd\mu_g}{(\varepsilon^2+\psi)^2}\\
 &=2\pi\int_0^{r_0}\frac{r}{(\varepsilon^2+r^2)^2}\dd r
   +O\left(\int_0^{r_0}\frac{r^3}{(\varepsilon^2+r^2)^2}\dd r\right)+O(1)\\
 &=\pi\left(\varepsilon^{-2}-(\varepsilon^2+r_0^2)^{-1}\right)
   +O\left(\log\frac{\varepsilon^2+r_0^2}{\varepsilon^2}\right)+O(1)\\
 &=\pi\varepsilon^{-2}+O(\log(\varepsilon^{-1})).
\end{align*}
Together with \eqref{eq:bubble-mean-bound}, this yields
\begin{equation}\label{eq:bubble-log-integral}
 \log\left(\int_\Sigma e^{\phi_\varepsilon}\dd\mu_g\right)
 =2\log(\varepsilon^{-1})+O(1).
\end{equation}
Since $|\nabla_gr|=1$ on $B_g(x_0,r_0)\setminus\{x_0\}$ and $\psi\ge r_0^2$ outside $B_g(x_0,r_0)$, we obtain
\begin{align*}
 \int_\Sigma|\nabla_g\phi_\varepsilon|^2\dd\mu_g
 &=\int_\Sigma\frac{4|\nabla_g\psi|^2}{(\varepsilon^2+\psi)^2}\dd\mu_g\\
 &=16\int_{B_g(x_0,r_0)}
       \frac{r^2}{(\varepsilon^2+r^2)^2}\dd\mu_g
   +4\int_{\Sigma\setminus B_g(x_0,r_0)}
       \frac{|\nabla_g\psi|^2}{(\varepsilon^2+\psi)^2}\dd\mu_g\\
 &=32\pi\int_0^{r_0}\frac{r^3}{(\varepsilon^2+r^2)^2}\dd r
   +O\left(\int_0^{r_0}\frac{r^5}{(\varepsilon^2+r^2)^2}\dd r\right)+O(1)\\
 &=32\pi\int_0^{r_0}\frac{r^3}{(\varepsilon^2+r^2)^2}\dd r
   +O\left(\int_0^{r_0}r\dd r\right)+O(1)\\
 &=16\pi\left(\log\frac{\varepsilon^2+r_0^2}{\varepsilon^2}
       -\frac{r_0^2}{\varepsilon^2+r_0^2}\right)+O(1)\\
 &=32\pi\log(\varepsilon^{-1})+O(1).
\end{align*}
Combining this with \eqref{eq:bubble-initial-energy}, \eqref{eq:bubble-log-integral} and $M>8\pi$, we obtain
\[
 \mathcal E_g(u_{0,\varepsilon},v_{0,\varepsilon})
 =2(8\pi-M)\log(\varepsilon^{-1})+O(1)
 \rightarrow-\infty
 \quad\text{as }\varepsilon\searrow0.
\]
Choosing $\varepsilon$ sufficiently small and setting $u_0=u_{0,\varepsilon}$ and $w_0=w_{0,\varepsilon}$ completes the proof.
\end{proof}

We can now prove Theorem~\ref{thm:supercritical}.

\begin{proof}[Proof of Theorem~\ref{thm:supercritical}]
Let $\mathcal E_g$ be the energy functional defined in \eqref{eq:energy}.
Recall that $\mathcal S_M$ is the set of smooth stationary solutions $(u_*,v_*,w_*)$ of \eqref{eq:system} satisfying $u_*,w_*\ge0$, $\int_\Sigma u_*\dd\mu_g=M$, and $\int_\Sigma v_*\dd\mu_g=0$.
Set
\[
 \beta_M:=\inf_{(u_*,v_*,w_*)\in\mathcal S_M}\mathcal E_g(u_*,v_*),
\]
which is finite by Proposition~\ref{prop:stationary-energy-lower}.
Apply Lemma~\ref{lem:supercritical-low-energy} with $R=\beta_M-1$ to choose smooth positive initial data $u_0,w_0$ with $\int_\Sigma u_0\dd\mu_g=M$ such that
\[
 \mathcal E_g(0)<\beta_M-1.
\]
Proposition~\ref{prop:global-existence} provides a global classical solution $(u,v,w)$.
Suppose that this solution satisfies $\sup_{t\ge0}\norm{u(t)}_{L^\infty}<\infty$.
By Lemma~\ref{lem:bounded-stationary-limit}, there exist a sequence $(t_j)_{j\ge1}$ with $t_j\to\infty$ as $j\to\infty$ and $(u_*,v_*,w_*)\in\mathcal S_M$ such that
\begin{equation*}
 \beta_M\le\mathcal E_g(u_*,v_*)
 =\lim_{j\to\infty}\mathcal E_g(t_j)
 \le\mathcal E_g(0)<\beta_M-1,
\end{equation*}
which is a contradiction.
Hence $\sup_{t\ge0}\norm{u(t)}_{L^\infty}=\infty$.
This completes the proof.
\end{proof}

\section{Unbounded solutions of the planar problem for supercritical mass}
\label{sec:planar-neumann}
In this section, we prove Theorem~\ref{thm:planar-infinite-time-blow-up}.
For the planar problem, the energy functional corresponding to \eqref{eq:energy} is
\begin{equation*}
\mathcal E(u,v):=\int_\Omega(L(u)-uv)\dd x
+\frac12\int_\Omega|\nabla v|^2\dd x,
\end{equation*}
where $L$ is defined in \eqref{eq:entropy-density}.
The following proposition shows that, for each $M>0$, the energy is bounded below uniformly over all stationary solutions of mass $M$ on $\Omega$.
\begin{proposition}\label{prop:planar-stationary-energy}
Let $M>0$ and let $\mathcal S_M$ denote the set of smooth steady states $(u,v,w)$ of \eqref{eq:planar-evolution} satisfying \eqref{eq:stationary-triple} and \eqref{eq:stationary-mean-field} on $\Omega$.
Then
\begin{equation*}
\inf_{(u,v,w)\in\mathcal S_M}
 \mathcal E(u,v)>-\infty.
\end{equation*}
\end{proposition}
To prove this Proposition, let $\Sigma$ be the Schottky double of $\overline\Omega$ \cite[Section~2.1]{2012-IUMJ-GustafssonSebbar}, a compact connected Riemann surface without boundary.
More explicitly, $\Sigma:=\Omega\cup\partial\Omega\cup\widetilde\Omega$, where $\widetilde\Omega$ is an identical copy of $\Omega$.
Let $\tau:\Sigma\to \Sigma$ be the antiholomorphic involution that sends each point of $\Omega$ to its corresponding point in $\widetilde\Omega$ and fixes $\partial\Omega$ pointwise.

\begin{lemma}\label{lem:planar-double-metric}
There is a smooth positive function $\rho$ on $\overline\Omega$ such that on $\Omega$,
\begin{equation}\label{eq:planar-conformal-factor}
g_{\mathrm{E}}=\rho g,
\qquad \Delta=\rho^{-1}\Delta_g.
\end{equation}
The even extension of $\rho$ to $\Sigma$ belongs to $W^{1,\infty}(\Sigma)$, has $\inf_\Sigma\rho>0$, and satisfies
\begin{equation}\label{eq:planar-factor-integral}
\int_\Sigma\rho\dd\mu_g=2|\Omega|.
\end{equation}
\end{lemma}
\begin{proof}
By \cite[Lemma~2.3.3]{2006-Jost}, $\Sigma$ admits a smooth conformal metric $g_0$ \cite[Definition~2.3.1]{2006-Jost}.
Then, $g:=(g_0+\tau^*g_0)/2$
is conformal and $\tau$ is an isometry.
The local formulas in \cite[Definitions~2.3.1 and~2.3.3]{2006-Jost} give \eqref{eq:planar-conformal-factor}, from which \eqref{eq:planar-factor-integral} follows directly.
\end{proof}
\begin{lemma}\label{lem:planar-correction}
Fix $g$ and $\rho$ as in Lemma~\ref{lem:planar-double-metric}.
For $(u,v,w)\in\mathcal S_M$, set $Z:=\int_\Omega e^v\dd x$, extend $v$ evenly to $\Sigma$, and define
\[
\eta:=v-\log(2Z).
\]
Then $\eta\in C^2(\Sigma)$ satisfies
\begin{equation}\label{eq:planar-transformed}
-\Delta_g\eta=2M\left(\rho e^\eta-\frac{\rho}{2A}\right),
\qquad \int_\Sigma\rho e^\eta\dd\mu_g=1.
\end{equation}
Moreover, there is $C=C(\Omega,M)>0$ such that the spatial average of the even extension of $v$ on $\Sigma$ satisfies
\begin{equation}\label{eq:planar-mean-bound}
|\overline v|\le C\qquad\text{for all }(u,v,w)\in\mathcal S_M.
\end{equation}
\end{lemma}

\begin{proof}
The Neumann condition ensures that the even extension of $v$ belongs to $C^2(\Sigma)$.
Combining \eqref{eq:stationary-mean-field} and \eqref{eq:planar-conformal-factor}, we obtain
\[
-\Delta_g v
=2M\left(\rho e^\eta-\frac{\rho}{2A}\right)\quad \text{on } \Sigma,
\]
and \eqref{eq:planar-factor-integral} shows that
\[
\int_\Sigma\rho e^\eta\dd\mu_g
=2\int_\Omega\frac{e^v}{2Z}\dd x=1,
\]
which yields \eqref{eq:planar-transformed}.
Let $G_g$ be the Green kernel of $-\Delta_g$.
The logarithmic bound in \cite[Theorem~4.13]{1982-Aubin} implies
\[
\sup_{y\in \Sigma}\int_\Sigma|G_g(x,y)|\dd\mu_g(x)\le C_1.
\]
Then, together with \eqref{eq:planar-factor-integral} and \eqref{eq:planar-transformed}, we obtain
\[
\norm{v-\overline v}_{L^1(\Sigma)}
\le 2MC_1\int_\Sigma\rho\left|e^\eta-\frac{1}{2A}\right|\dd\mu_g
\le 4C_1M.
\]
Finally, $\int_\Sigma\rho v\dd\mu_g=2\int_\Omega v\dd x=0$ and \eqref{eq:planar-factor-integral} imply
\[
2A|\overline v|
=\left|\int_\Sigma\rho(v-\overline v)\dd\mu_g\right|
\le4C_1M\norm\rho_{L^\infty(\Sigma)},
\]
which proves \eqref{eq:planar-mean-bound}.
\end{proof}

We now prove Proposition~\ref{prop:planar-stationary-energy}.

\begin{proof}[Proof of Proposition~\ref{prop:planar-stationary-energy}]
Fix $(u,v,w)\in\mathcal S_M$ and set $Z:=\int_\Omega e^v\dd x$.
Testing the first equation in \eqref{eq:stationary-mean-field} by $v$ gives
\[
\int_\Omega|\nabla v|^2\dd x
=M\int_\Omega v\frac{e^v}{Z}\dd x.
\]
Since $\log u=\log M+v-\log Z$, we obtain
\begin{equation}\label{eq:planar-energy-identity}
\mathcal E(u,v)=M\log M-M+A-M\log Z
+\frac M2\int_\Omega v\frac{e^v}{Z}\dd x.
\end{equation}
With $\eta$ as in Lemma~\ref{lem:planar-correction}, reflection symmetry gives $\int_\Sigma v\rho e^\eta\dd\mu_g=\int_\Omega ve^v/Z\dd x$.
Using $v=\eta+\log(2Z)$ and $\overline\eta=\overline v-\log(2Z)$ in \eqref{eq:planar-energy-identity}, we find
\begin{equation*}
\mathcal E(u,v)=M\log M-M+A
+\frac M2\left(\int_\Sigma\eta\rho e^\eta\dd\mu_g+\overline\eta\right)
+M\log2-\frac M2\overline v.
\end{equation*}
Applying Lemma~\ref{lem:weighted-energy} with $\lambda=2M$ and using Lemma~\ref{lem:planar-correction}, we complete the proof.
\end{proof}

\begin{proof}[Proof of Theorem~\ref{thm:planar-infinite-time-blow-up}]
We claim that for every $M>4\pi$ and $R\in\mathbb R$,
there exist strictly positive $u_0,w_0\in C^\infty(\overline\Omega)$ satisfying $\int_\Omega u_0\dd x=M$ and
\[
 \mathcal E(u_0,v_0)<R.
\]
We follow the construction in Lemma~\ref{lem:supercritical-low-energy}, with concentration at a boundary point.
Fix $M>4\pi$ and $x_0\in\partial\Omega$, and choose boundary normal coordinates $(s,y)$ centered at $x_0$.
Let $r:=(s^2+y^2)^{1/2}$ and let $r_0>0$ so that $H_{2r_0}:=\{r<2r_0,\ y\ge0\}$ lies in the chart.
Choose a smooth radial cutoff $\chi$ equal to one for $r\le r_0$ and zero for $r\ge2r_0$, and set
\[
 \psi:=\chi r^2+(1-\chi)r_0^2
\]
in the chart, extending it by $r_0^2$ outside.
Then $\partial_\nu\psi=0$, $\psi=r^2$ on $H_{r_0}$, and $\psi\ge r_0^2$ outside $H_{r_0}$.
For $0<\varepsilon<\min\{1,r_0\}$, define $\phi_\varepsilon$ and the initial data $(u_{0,\varepsilon},w_{0,\varepsilon})$ as in the proof of Lemma~\ref{lem:supercritical-low-energy}, using the above $\psi$ on $\overline\Omega$.
The estimates \eqref{eq:bubble-mean-bound} and \eqref{eq:bubble-log-integral} hold with $\Sigma$ replaced by $\Omega$.
Using polar coordinates with $0\le\theta\le\pi$ on $H_{r_0}$, we obtain
\begin{align*}
 \int_\Omega|\nabla\phi_\varepsilon|^2\dd x
 &=16\pi\int_0^{r_0}\frac{r^3}{(\varepsilon^2+r^2)^2}\dd r
   +O\left(\int_0^{r_0}\frac{r^4}{(\varepsilon^2+r^2)^2}\dd r\right)+O(1)\\
 &=8\pi\left(\log\frac{\varepsilon^2+r_0^2}{\varepsilon^2}
       -\frac{r_0^2}{\varepsilon^2+r_0^2}\right)+O(1)\\
 &=16\pi\log(\varepsilon^{-1})+O(1).
\end{align*}
Substituting this estimate into the Euclidean version of \eqref{eq:bubble-initial-energy} gives $\mathcal E(u_{0,\varepsilon},v_{0,\varepsilon})=2(4\pi-M)\log(\varepsilon^{-1})+O(1)\to-\infty$ as $\varepsilon\searrow0$.
Taking $\varepsilon$ sufficiently small proves the claim.

The proof of Proposition~\ref{prop:global-existence} gives global classical solvability on $\Omega$, using Neumann elliptic and heat semigroup estimates.
The arguments of Lemmas~\ref{lem:energy-monotonicity} and~\ref{lem:bounded-stationary-limit} give energy monotonicity and stationary limits of bounded solutions.
Then, together with Proposition~\ref{prop:planar-stationary-energy}, we complete the proof.
\end{proof}

\section{Uniform boundedness at critical mass on spheres}
\label{sec:critical-sphere}

In this section, we prove Theorem~\ref{thm:critical-sphere}.
Let $\Sigma=\mathbb S_R$ with the canonical spherical metric, $A=4\pi R^2$, and $M=8\pi$.

\begin{lemma}\label{lem:concentration-rearrangement}
For a nonnegative $f\in C(\Sigma)$, define its mass concentration function \cite[Definition~2.3]{2015-ASNSP-Langford} by
\begin{equation}\label{eq:spherical-concentration}
 \mathcal C_f(s):=
 \sup\left\{\int_E f\dd\mu_g:
 E\subset\Sigma\text{ measurable},\ \mu_g(E)=As\right\},
 \qquad 0\le s\le1.
\end{equation}
Following \cite[Definition~2.2]{2015-ASNSP-Langford}, the decreasing rearrangement $f^*:(0,A)\to[0,\infty)$ is defined by
\begin{equation}\label{eq:decreasing-rearrangement-definition}
 f^*(r):=\inf\left\{\lambda\ge0:
 \mu_g\bigl(\{x\in\Sigma:f(x)>\lambda\}\bigr)\le r\right\},
 \qquad 0<r<A.
\end{equation}
By \cite[Proposition~2.6]{2015-ASNSP-Langford}, for $0<s<1$,
\begin{equation*}
 \mathcal C_f(s)=\int_0^{As}f^*(r)\dd r.
\end{equation*}
Consequently,
\begin{equation}\label{eq:concentration-derivative}
 (\mathcal C_f)'(s)=Af^*(As)\quad\text{for all }0<s<1.
\end{equation}
\end{lemma}

We next derive a differential inequality for the mass concentration function of $u$.

\begin{lemma}\label{lem:spherical-concentration}
Let $(u,v,w)$ be a global classical solution of \eqref{eq:system}--\eqref{eq:normalization} with mass $M=8\pi$, and set
\[
 U(s,t):=\mathcal C_{u(t)}(s),\qquad
 W(s,t):=\mathcal C_{w(t)}(s).
\]
Then $U(0,t)=0$, $U(1,t)=M$, $\partial_s U\ge0$, and
\begin{equation}\label{eq:spherical-concentration-pde}
 \partial_t U\le\frac{s(1-s)}{R^2}\partial_s^2 U
       +\left(\frac{W(s,t)}{A}-s\overline w\right)\partial_s U
 \quad\text{for a.e. }(s,t)\in(0,1)\times(0,\infty).
\end{equation}
For every $0<\tau<T<\infty$, $U$ is Lipschitz on $[0,1]\times[\tau,T]$ and $\partial_s^2 U\in L^1((0,1)\times(\tau,T))$.
\end{lemma}

\begin{proof}
For any nonnegative $f\in C(\Sigma)$, the definition \eqref{eq:spherical-concentration} gives
\begin{equation}\label{eq:rearrangement-quantile-area}
 \mu_g\bigl(\{f>f^*(As)\}\bigr)\le As
 \le\mu_g\bigl(\{f\ge f^*(As)\}\bigr)
 \quad\text{for all }0<s<1,
\end{equation}
and
\begin{equation}\label{eq:rearrangement-level-lower}
 f^*(r)>t\quad\Longleftrightarrow\quad r<\mu_g(\{f>t\}),
 \qquad t\ge0,\quad 0<r<A.
\end{equation}
Fix $0<\tau<T<\infty$.
For $t,r\in[\tau,T]$ and any measurable set $E$ with $\mu_g(E)=As$, the definition \eqref{eq:spherical-concentration} gives
\begin{align*}
 \int_Eu(t)\dd\mu_g
 &=\int_Eu(r)\dd\mu_g+\int_E\bigl(u(t)-u(r)\bigr)\dd\mu_g\\
 &\le U(s,r)+As\norm{u(t)-u(r)}_{L^\infty}
 \quad\text{for all }0\le s\le1,\ t,r\in[\tau,T].
\end{align*}
Taking the supremum over $E$ and then interchanging $t$ and $r$ yields, for all $0\le s\le1,\ t,r\in[\tau,T]$,
\begin{align}
 & |U(s,t)-U(s,r)|
 \le As\norm{u(t)-u(r)}_{L^\infty}
 \le A\norm{u(t)-u(r)}_{L^\infty}.
\label{eq:concentration-time-estimate}
\end{align}
By \eqref{eq:decreasing-rearrangement-definition} and \eqref{eq:concentration-derivative},
\begin{equation}\label{eq:concentration-spatial-estimate}
\partial_s U(s,t)=A\,[u(t)]^*(As)
 \le A\norm{u(t)}_{L^\infty}\quad\text{for all }0<s<1,\ t\in[\tau,T].
\end{equation}
Together with the regularity of $u$, \eqref{eq:concentration-time-estimate} and \eqref{eq:concentration-spatial-estimate} show that $U$ is Lipschitz on $[0,1]\times[\tau,T]$.

Define
\begin{equation}\label{eq:concentration-threshold-definition}
 \lambda(s,t):=\frac{\partial_s U(s,t)}{A}.
\end{equation}

Suppose first that $\lambda(s,t)$ is a regular value of $u(t)$, and write $\lambda=\lambda(s,t)$ below.
By \eqref{eq:rearrangement-quantile-area} and \eqref{eq:concentration-spatial-estimate},
\[
 \mu_g(\{u(t)>\lambda\})\le As
 \le\mu_g(\{u(t)\ge\lambda\})
 \quad\text{for all }0<s<1,\ t\in[\tau,T].
\]
Hence, 
$E:=\{u(t)>\lambda\}$ satisfies $\mu_g(E)=As$.
Hence $E$ attains the supremum in \eqref{eq:spherical-concentration}, that is, 
\[
 U(s,t)=\int_Eu(t)\dd\mu_g.
\]
With $E$ fixed, $U(s,r)-\int_Eu(x,r)\dd\mu_g$ is nonnegative for $r\in[\tau,T]$ and vanishes at $r=t$.
Thus, by the Lipschitz continuity of $U$, for almost every $(s,t)\in(0,1)\times(\tau,T)$, we obtain
\begin{equation}\label{eq:concentration-time-derivative}
 \partial_t U(s,t)=\int_E\partial_t u(x,t)\dd\mu_g.
\end{equation}
Let $\nu$ be the outward unit normal to $\partial E$, and let $\mathrm{d}\ell_g$ denote the measure on $\partial E$.
Since $\nu=-\nabla_gu/|\nabla_gu|$ and $\partial_\nu u=-|\nabla_gu|$ on $\partial E$, \eqref{eq:system} and the above identity give
\begin{align}
 \partial_t U(s,t)
 &=\int_E\bigl(\Delta_gu-\operatorname{div}_g(u\nabla_gv)\bigr)\dd\mu_g
 \notag\\
 &=\int_{\partial E}\partial_\nu u\mathrm{d}\ell_g
   -\int_{\partial E}u\partial_\nu v\mathrm{d}\ell_g\notag\\
 &=-\int_{\partial E}|\nabla_gu|\mathrm{d}\ell_g
   -\lambda\int_{\partial E}\partial_\nu v\mathrm{d}\ell_g\notag\\
 &=-\int_{\partial E}|\nabla_gu|\mathrm{d}\ell_g
   -\lambda\int_E\Delta_gv\dd\mu_g\notag\\
 &=-\int_{\partial E}|\nabla_gu|\mathrm{d}\ell_g
   +\lambda\int_E\left(w-\overline w\right)\dd\mu_g.
 \label{eq:concentration-level-identity}
\end{align}
Fix $t\in[\tau,T]$.
We apply the coarea formula \cite[Corollary~I.3.1]{2001-Chavel} to see that for sufficiently small $h>0$,
\[
 \mu_g(\{u(t)>\lambda\})-\mu_g(\{u(t)>\lambda+h\})
 =\int_\lambda^{\lambda+h}
   \left(\int_{\{u(t)=\eta\}}|\nabla_gu|^{-1}\mathrm{d}\ell_g\right)\dd\eta,
\]
which yields, for $\eta$ near $\lambda$,
\[
 \frac{\dd}{\dd\eta}\mu_g(\{u(t)>\eta\})
 =-\int_{\{u(t)=\eta\}}|\nabla_gu|^{-1}\mathrm{d}\ell_g<0.
\]
Differentiating $\mu_g(\{u(t)>\lambda(s,t)\})=As$ with respect to $s$ and using $\partial_s U(s,t)=A\lambda(s,t)$, we obtain
\[
 \partial_s^2 U(s,t)=A\,\partial_s\lambda(s,t)
 =-\frac{A^2}{\int_{\partial E}|\nabla_gu|^{-1}\mathrm{d}\ell_g}.
\]
The Hölder inequality and spherical isoperimetric inequality \cite[(4.8)]{1978-BAMS-Osserman} yield
\begin{align}
 &\left(\int_{\partial E}|\nabla_gu|\mathrm{d}\ell_g\right)
 \left(\int_{\partial E}|\nabla_gu|^{-1}\mathrm{d}\ell_g\right)\notag\\
 &\qquad\ge\left(\int_{\partial E}\mathrm{d}\ell_g\right)^2
 \ge\frac{\mu_g(E)\bigl(A-\mu_g(E)\bigr)}{R^2}
 =\frac{A^2}{R^2}s(1-s).
 \label{eq:concentration-isoperimetric-estimate}
\end{align}
Combining these estimates with $\int_Ew\dd\mu_g\le W(s,t)$ in \eqref{eq:concentration-level-identity} proves \eqref{eq:spherical-concentration-pde} at almost every $(s,t)$ for which $\lambda$ is a regular value of $u(t)$.

To treat critical values, we first establish the local Lipschitz continuity of $\partial_s U$ and the integrability of $\partial_s^2 U$.
Fix $t\in[\tau,T]$, $\delta\in(0,\frac12)$, $\delta\le s_1<s_2\le1-\delta$, and $\eta$ with $\partial_s U(s_2,t)/A<\eta<\partial_s U(s_1,t)/A$.
Then \eqref{eq:concentration-derivative} and \eqref{eq:rearrangement-level-lower} imply
\[
 A\delta\le As_1<\mu_g(\{u(t)>\eta\})\le As_2\le A(1-\delta).
\]
Then, \eqref{eq:concentration-isoperimetric-estimate} gives
\begin{equation*}
 \int_{\{u(t)=\eta\}}\mathrm{d}\ell_g
 \ge C_1:=\frac{A\sqrt{\delta(1-\delta)}}{R}>0.
\end{equation*}
Combining this with the coarea formula \cite[Corollary~I.3.1]{2001-Chavel}, \eqref{eq:rearrangement-quantile-area} and \eqref{eq:concentration-spatial-estimate}, we obtain, for all $\delta\le s_1<s_2\le1-\delta$ and $t\in[\tau,T]$,
\begin{align}
 \frac{C_1}{A}\bigl(\partial_s U(s_1,t)-\partial_s U(s_2,t)\bigr)
 &\le\int_{\partial_s U(s_2,t)/A}^{\partial_s U(s_1,t)/A}
       \left(\int_{\{u(t)=\eta\}}\mathrm{d}\ell_g\right)\dd\eta\notag\\
 &=\int_{\{\partial_s U(s_2,t)/A<u(t)<\partial_s U(s_1,t)/A\}}
       |\nabla_gu(t)|\dd\mu_g\notag\\
 &\le\norm{\nabla_gu(t)}_{L^\infty}
       \mu_g\bigl(\{\partial_s U(s_2,t)/A<u(t)<\partial_s U(s_1,t)/A\}\bigr)\notag\\
 &\le A\norm{\nabla_gu(t)}_{L^\infty}(s_2-s_1).
 \label{eq:concentration-gradient-lipschitz}
\end{align}
Hence $\partial_s U$ is locally Lipschitz in $s$, uniformly for $t\in[\tau,T]$.
Moreover, by \eqref{eq:concentration-spatial-estimate}, $\partial_s U(\cdot,t)$ is nonincreasing, so $\partial_s^2 U(\cdot,t)\le0$ almost everywhere.
Therefore,
\begin{equation*}
 \int_a^b|\partial_s^2 U(s,t)|\dd s
 =-\int_a^b \partial_s^2 U(s,t)\dd s
 =\partial_s U(a,t)-\partial_s U(b,t),\qquad 0<a<b<1.
\end{equation*}
Letting $a\searrow0$ and $b\nearrow1$, and using \eqref{eq:concentration-spatial-estimate}, we obtain $\partial_s^2 U\in L^1((0,1)\times(\tau,T))$.
Fix $t\in[\tau,T]$.
With $\lambda$ as defined in \eqref{eq:concentration-threshold-definition}, let 
\[
 K_t:=\{u(x,t):x\in\Sigma,\ \nabla_gu(x,t)=0\},
\]
\[
 S_t:=\{s\in(0,1):\lambda(s,t)\in K_t\}.
\]
Since $u(t)\in C^2(\Sigma)$, Sard's theorem \cite[Proposition~I.3.8]{2001-Chavel}, applied in local coordinates, gives $|K_t|=0$.
By \eqref{eq:concentration-gradient-lipschitz}, $\lambda(\cdot,t)$ is Lipschitz on $[a,b]$ for any $0<a<b<1$.
We may therefore apply the one-dimensional area formula \cite[Theorem~IV.2.4]{2001-Chavel} to $S_t\cap(a,b)$ to obtain
\[
 \int_{S_t\cap(a,b)}|\partial_s\lambda(s,t)|\dd s
 =\int_{\mathbb R}
   \#\{s\in S_t\cap(a,b):\lambda(s,t)=y\}\dd y=0.
\]
Thus $\partial_s\lambda(\cdot,t)=0$ almost everywhere on $S_t\cap(a,b)$,
which implies $\partial_s^2 U=A\partial_s\lambda=0$ almost everywhere on $S_t$.
For $s\in S_t$, set $\lambda:=\lambda(s,t)$.
Choose a measurable set $E\subset\Sigma$ such that
\[
 \{u(t)>\lambda\}\subset E\subset\{u(t)\ge\lambda\},\qquad
 \mu_g(E)=As.
\]
Then, $E$ attains the supremum $U(s,t)$.
On $\{u=\lambda\}$, $\nabla_gu=0$ and $\Delta_gu=0$ almost everywhere.
Choose regular values $\lambda_n\searrow\lambda$ of $u(t)$.
The divergence theorem and dominated convergence give
\[
 \int_E\Delta_gu\dd\mu_g
 =\lim_{n\to\infty}\int_{\{u>\lambda_n\}}\Delta_gu\dd\mu_g
 =-\lim_{n\to\infty}\int_{\{u=\lambda_n\}}|\nabla_gu|\mathrm{d}\ell_g
 \le0.
\]
Moreover, integrating the divergence of $(u-\lambda)_+\nabla_gv$ over $\Sigma$ yields
\[
 \int_E\operatorname{div}_g((u-\lambda)\nabla_gv)\dd\mu_g=0.
\]
Thus, \eqref{eq:system} and \eqref{eq:concentration-time-derivative} give
\begin{align*}
 \partial_t U(s,t)
 &=\int_E\Delta_gu\dd\mu_g
   -\int_E\operatorname{div}_g(u\nabla_gv)\dd\mu_g\\
 &\le-\lambda\int_E\Delta_gv\dd\mu_g
 =\lambda\int_E(w-\overline w)\dd\mu_g\\
 &\le\lambda\bigl(W(s,t)-As\,\overline w\bigr).
\end{align*}
Together with $\partial_s^2 U=0$ almost everywhere on $S_t$, this proves \eqref{eq:spherical-concentration-pde} when $\lambda$ is critical.
\end{proof}

We now prove Theorem~\ref{thm:critical-sphere} using a comparison function related to that in \cite[Lemma~4.5]{2025-JEE-XiangYang}.

\begin{proof}[Proof of Theorem~\ref{thm:critical-sphere}]

By Proposition~\ref{prop:global-existence} and the strong maximum principle
\cite[Proposition~2.3]{2006-Mueller}, $u(\cdot,1)>0$ on $\Sigma$.
Taking $t=1$ as the new time origin, we may thus assume $\min_\Sigma u_0>0$.
Let $U,W$ be defined as in Lemma~\ref{lem:spherical-concentration}.
Since $W(s,0)\le\min\{A\norm{w_0}_{L^\infty}s,A\overline{w}_0\}$ and $0\le \overline{w}_0\le\norm{w_0}_{L^\infty}$, we have
\begin{align}
 W(s,0)-As\overline{w}_0
 &\le A\norm{w_0}_{L^\infty}\min\{s,1-s\}\notag\\
 &\le 2A\norm{w_0}_{L^\infty}s(1-s)\quad\text{for all }0\le s\le1.
 \label{eq:critical-initial-signal}
\end{align}
By \eqref{eq:w-representation}, we have, for all $0\le s\le1$ and $t\ge0$,
\[
\begin{aligned}
 & W(s,t)\le e^{-t}W(s,0)+\int_0^t e^{-(t-r)}U(s,r)\dd r,
 \qquad
 \overline w=e^{-t}\overline{w}_0+\frac MA(1-e^{-t}).
\end{aligned}
\]
Using $1-e^{-t}=\int_0^t e^{-(t-r)}\dd r$ and \eqref{eq:critical-initial-signal}, we obtain, for all $0\le s\le1$ and $t\ge0$,
\begin{align}
 W(s,t)-As\overline w
 &\le e^{-t}W(s,0)+\int_0^t e^{-(t-r)}U(s,r)\dd r
      -As e^{-t}\overline{w}_0-Ms(1-e^{-t})\notag\\
 &=e^{-t}\bigl(W(s,0)-As\overline{w}_0\bigr)
   +\int_0^t e^{-(t-r)}\bigl(U(s,r)-Ms\bigr)\dd r\notag\\
 &\le 2A\norm{w_0}_{L^\infty}e^{-t}s(1-s)
   +\int_0^t e^{-(t-r)}\bigl(U(s,r)-Ms\bigr)\dd r.
 \label{eq:critical-signal-bound}
\end{align}
Define
\[
 b(t):=\frac{\min_\Sigma u_0}{2\max_\Sigma u_0}
 \exp\left(-2\norm{w_0}_{L^\infty}(1-e^{-t})\right).
\]
Then
\begin{equation}\label{eq:critical-parameter-properties}
 \frac{\min_\Sigma u_0}{2\max_\Sigma u_0}e^{-2\norm{w_0}_{L^\infty}}\le b(t)\le\frac12,\qquad
 b'(t)=-2\norm{w_0}_{L^\infty}e^{-t}b(t).
\end{equation}
Set
\[
 \overline{U}(s,t):=\frac{Ms}{b(t)+(1-b(t))s}.
\]
Differentiating and using \eqref{eq:critical-parameter-properties}, we obtain
\begin{align}
 \partial_s \overline{U}(s,t)&=\frac{Mb(t)}{[b(t)+(1-b(t))s]^2},\qquad
 \partial_s^2 \overline{U}(s,t)=-\frac{2Mb(t)(1-b(t))}{[b(t)+(1-b(t))s]^3},
 \label{eq:critical-profile-s}\\
 \partial_t \overline{U}(s,t)&=-\frac{Ms(1-s)b'(t)}{[b(t)+(1-b(t))s]^2}
 =2\norm{w_0}_{L^\infty}e^{-t}s(1-s)\partial_s \overline{U}(s,t).
 \label{eq:critical-profile-t}
\end{align}
Then, with $M=8\pi$ and $A=4\pi R^2$, we obtain
\begin{equation}\label{eq:critical-profile-identity}
 \frac{s(1-s)}{R^2}\partial_s^2 \overline{U}+\frac{\overline{U}-Ms}{A}\partial_s \overline{U}=0.
\end{equation}
For $0<s<1$, the bounds $U(s,0)\le A(\max_\Sigma u_0)s$ and $M-U(s,0)\ge A(\min_\Sigma u_0)(1-s)$ imply
\[
 \frac{U(s,0)}{M-U(s,0)}
 \le\frac{\max_\Sigma u_0}{\min_\Sigma u_0}\frac{s}{1-s}
 \le\frac1{b(0)}\frac{s}{1-s}
 =\frac{\overline{U}(s,0)}{M-\overline{U}(s,0)}
 \quad\text{for all }0<s<1,
\]
which yields
\begin{equation}\label{eq:critical-initial-order}
 U(s,0)\le \overline{U}(s,0)\quad\text{for all }0\le s\le1.
\end{equation}
Since $\overline{U}\ge Ms$ and $\partial_s \overline{U},\partial_t \overline{U}\ge0$ by \eqref{eq:critical-profile-t} and \eqref{eq:critical-profile-identity}, we have, for all $0\le s\le1$ and $t\ge0$,
\begin{align}
 \partial_t \overline{U}(s,t)-\frac{s(1-s)}{R^2}\partial_s^2 \overline{U}(s,t)
 &=\frac{\partial_s \overline{U}(s,t)}{A}
 \bigl(\overline{U}(s,t)-Ms+2A\norm{w_0}_{L^\infty}e^{-t}s(1-s)\bigr)\notag\\
 &\ge\frac{\partial_s \overline{U}(s,t)}{A}\biggl[
 2A\norm{w_0}_{L^\infty}e^{-t}s(1-s)\notag\\
 &\qquad+\int_0^t e^{-(t-r)}\bigl(\overline{U}(s,r)-Ms\bigr)\dd r\biggr].
 \label{eq:critical-memory-supersolution}
\end{align}
Let $D:=U-\overline{U}$.
Combining \eqref{eq:spherical-concentration-pde}, \eqref{eq:critical-signal-bound} and \eqref{eq:critical-memory-supersolution} gives, for almost every $(s,t)\in(0,1)\times(0,\infty)$,
\begin{align}
 & \partial_t D\le\frac{s(1-s)}{R^2}\partial_s^2 D
 +\left(\frac{W}{A}-s\overline w\right)\partial_s D
 +\frac{\partial_s \overline{U}}{A}\int_0^t e^{-(t-r)}D(s,r)\dd r.
\label{eq:critical-difference}
\end{align}
By \eqref{eq:critical-initial-order}, $D(s,0)\le0$, and $D(0,t)=D(1,t)=0$.
Fix $T>0$ and set $D_+:=\max\{D,0\}$.
Using the regularity of $U$ and the integrability of $\partial_s^2 U$ from Lemma~\ref{lem:spherical-concentration}, we test \eqref{eq:critical-difference} with $D_+$ and integrate by parts.
By \eqref{eq:critical-profile-s}, $\partial_s W\ge0$, and Young's inequality, we obtain, for almost every $t\in(0,T)$,
\begin{align*}
 \frac12\frac{\dd}{\dd t}\norm{D_+(\cdot,t)}_{L^2(0,1)}^2
 &\le\frac1{R^2}\int_0^1s(1-s)(\partial_s^2 D)D_+\dd s
   +\int_0^1\left(\frac WA-s\overline w\right)(\partial_s D)D_+\dd s\\
 &\quad+\frac1A\int_0^1\partial_s \overline{U}(s,t)D_+(s,t)
       \int_0^t e^{-(t-r)}D_+(s,r)\dd r\mathrm{d}s\\
 &=-\frac1{R^2}\int_0^1s(1-s)|\partial_sD_+|^2\dd s
   -\frac1{2R^2}\int_0^1(1-2s)\partial_s(D_+^2)\dd s\\
 &\quad+\frac12\int_0^1\left(\frac WA-s\overline w\right)\partial_s(D_+^2)\dd s\\
 &\quad+\frac1A\int_0^1\partial_s \overline{U}(s,t)D_+(s,t)
       \int_0^t e^{-(t-r)}D_+(s,r)\dd r\mathrm{d}s\\
 &=-\frac1{R^2}\int_0^1s(1-s)|\partial_sD_+|^2\dd s
   +\frac12\int_0^1\left(\overline w-\frac2{R^2}-\frac{\partial_s W}A\right)D_+^2\dd s\\
 &\quad+\frac1A\int_0^1\partial_s \overline{U}(s,t)D_+(s,t)
       \int_0^t e^{-(t-r)}D_+(s,r)\dd r\mathrm{d}s\\
 &\le\frac{\overline w}2\norm{D_+(\cdot,t)}_{L^2(0,1)}^2\\
 &\quad+\frac1{2A}\int_0^t e^{-(t-r)}
       \int_0^1 \partial_s \overline{U}(s,t)\bigl(D_+(s,t)^2+D_+(s,r)^2\bigr)\dd s\mathrm{d}r\\
 &\le\left(\frac{\overline w}2+\frac{M(1-e^{-t})}{2Ab(t)}\right)
       \norm{D_+(\cdot,t)}_{L^2(0,1)}^2\\
 &\quad+\frac{M}{2Ab(t)}\int_0^t e^{-(t-r)}
       \norm{D_+(\cdot,r)}_{L^2(0,1)}^2\dd r\\
 &\le\frac{C_1}{2}\norm{D_+(\cdot,t)}_{L^2(0,1)}^2
       +\frac{C_1}{2}\int_0^t\norm{D_+(\cdot,r)}_{L^2(0,1)}^2\dd r,
\end{align*}
where $C_1:=\max_{0\le t\le T}\overline w+M/(A\min_{0\le t\le T}b(t))>0$.
Since $D_+(\cdot,0)=0$, applying Gronwall's inequality yields $D_+\equiv0$ on $[0,1]\times[0,T]$.
Since $T$ was arbitrary,
\begin{equation*}
 U(s,t)\le \overline{U}(s,t)\quad\text{for all }0\le s\le1\text{ and }t\ge0.
\end{equation*}
Using this and the lower bound for $b(t)$ in \eqref{eq:critical-parameter-properties}, we obtain
\begin{align*}
 \norm{u(t)}_{L^\infty}
 &=\lim_{s\searrow0}\frac{U(s,t)}{As}
 \le\frac{M}{Ab(t)}\\
 &\le\frac{4\max_\Sigma u_0}{R^2\min_\Sigma u_0}
 e^{2\norm{w_0}_{L^\infty}}
 \quad\text{for all }t\ge0.
\end{align*}
The representation \eqref{eq:w-representation} then gives
\[
 \sup_{t\ge0}\norm{w(t)}_{L^\infty}
 \le\max\left\{\norm{w_0}_{L^\infty},
 \frac{4\max_\Sigma u_0}{R^2\min_\Sigma u_0}
 e^{2\norm{w_0}_{L^\infty}}\right\}.
\]
The uniform $W^{1,\infty}$ bound for $v$ follows from the elliptic estimate \cite[Theorem~A.3]{2019-Lee} and Sobolev embedding \cite[Theorem~2.7]{1999-Hebey}.
Finally, by continuity, the solution is bounded on $[0,1]$.
This completes the proof.
\end{proof}

\section{Unbounded solutions for critical mass}
\label{sec:critical-grow-up}

To prove Theorem~\ref{thm:critical-grow-up}, we first verify that $g$ extends smoothly to $\Sigma$ and then construct a comparison function.
\begin{lemma}\label{lem:critical-grow-up-metric}
Let $\Sigma$ and $g$ be as in Theorem~\ref{thm:critical-grow-up}.
Then $g$ extends to a smooth rotation-invariant metric on $\Sigma$, 
and
\begin{equation}\label{eq:area-formulas}
 \mathrm{d}\mu_g=2\,\mathrm{d}s\,\mathrm{d}\theta,\qquad
 \mu_g(\Sigma)=4\pi,\qquad
 \mu_g(C_s)=4\pi s,
\end{equation}
where
\[
 C_s:=
 \bigl\{(\sigma,\theta):0\le\sigma\le s,\ \theta\in\mathbb R/(2\pi\mathbb Z)\bigr\},
 \qquad 0<s<1.
\]
Moreover, for a rotation-invariant smooth function $f=f(s)$,
\begin{equation}\label{eq:radial-laplacian}
\Delta_g f=\partial_s(D\partial_s f).
\end{equation}
\end{lemma}

\begin{proof}
By \eqref{eq:critical-grow-up-coordinate}, $s\to0$ as $z\to0$
and $s\to1$ as $|z|\to\infty$.
In the coordinate $z$, the metric takes the form
\begin{equation}\label{eq:metric-complex-coordinate}
g=4D(s)\bigl[(\mathrm{d}\log|z|)^2+\mathrm{d}\theta^2\bigr]
=\frac{4D(s)}{|z|^2}|\mathrm{d}z|^2.
\end{equation}
Since $D(s)/s$ is analytic and equals $1$ at $0$, we have
\[
|z|^2=sA(s),\qquad
A(s):=2\exp\!\left(\int_{1/2}^{s}
\left[\frac1{D(\sigma)}-\frac1\sigma\right]\mathrm{d}\sigma\right),
\]
where $A$ is analytic near $0$ and $A(0)>0$.
The inverse function theorem implies that $s$ is analytic in $|z|^2$ near $z=0$.
By \eqref{eq:metric-complex-coordinate} and
\[
\lim_{s\to0}\frac{4D(s)}{sA(s)}
=\frac4{A(0)}>0,
\]
we see that $g$ extends smoothly and remains positive definite at $z=0$.
Since $D(1-s)=D(s)$, replacing $s,z$ by $1-s,1/z$ gives the same conclusion at $z=\infty$.
The independence of $D$ from $\theta$ implies that $g$ is rotation-invariant.
Moreover, $\det g=4$ gives \eqref{eq:area-formulas}.
For $f=f(s)$,
\[
\Delta_g f
=\frac1{\sqrt{\det g}}\partial_s\!\left(\sqrt{\det g}\,D \partial_s f\right)
=\partial_s(D\partial_s f),
\]
which yields \eqref{eq:radial-laplacian}.

\end{proof}

\begin{lemma}\label{lem:critical-grow-up-cumulative}
Let $(u,v,w)$ be a rotation-invariant global classical solution of \eqref{eq:system} on $(\Sigma,g)$ with
\begin{equation}\label{eq:initial-masses}
\int_\Sigma u_0\,\mathrm{d}\mu_g=\int_\Sigma w_0\,\mathrm{d}\mu_g=8\pi.
\end{equation}
Define
\[
\widehat{U}(s,t):=\frac12\int_0^s u(\sigma,t)\,\mathrm{d}\sigma,
\qquad
\widehat{W}(s,t):=\frac12\int_0^s w(\sigma,t)\,\mathrm{d}\sigma.
\]
Then
\begin{equation}\label{eq:cumulative-system}
\partial_t \widehat{U}=D \partial_s^2 \widehat{U}+2(\widehat{W}-s)\partial_s \widehat{U},\qquad \partial_t \widehat{W}=\widehat{U}-\widehat{W}
\quad\text{on }(0,1)\times(0,\infty),
\end{equation}
with $\widehat{U}(0,t)=\widehat{W}(0,t)=0$ and $\widehat{U}(1,t)=\widehat{W}(1,t)=1$ for all $t\ge0$.
Moreover,
\begin{equation}\label{eq:Q-representation}
\widehat{W}(s,t)=e^{-t}\widehat{W}(s,0)+\int_0^t e^{-(t-r)}\widehat{U}(s,r)\,\mathrm{d}r.
\end{equation}
\end{lemma}

\begin{proof}
Integration of \eqref{eq:system} and application of the divergence theorem, together with \eqref{eq:initial-masses}, yield
\[
\frac{\mathrm{d}}{\mathrm{d}t}\int_\Sigma u\,\mathrm{d}\mu_g=0,\qquad
\frac{\mathrm{d}}{\mathrm{d}t}\int_\Sigma w\,\mathrm{d}\mu_g
=8\pi-\int_\Sigma w\,\mathrm{d}\mu_g,
\]
and hence
\[
\int_\Sigma u(t)\,\mathrm{d}\mu_g=\int_\Sigma w(t)\,\mathrm{d}\mu_g=8\pi.
\]
Since $\mu_g(\Sigma)=4\pi$ by Lemma~\ref{lem:critical-grow-up-metric}, we obtain
\[
\overline w=\frac{1}{4\pi}\int_\Sigma w(t)\,\mathrm{d}\mu_g=2.
\]
The definitions of $\widehat{U}$ and $\widehat{W}$ give
\[
u=2\partial_s \widehat{U},\qquad w=2\partial_s \widehat{W},
\]
and
$$
\partial_t \widehat{W}=\frac12\int_0^s(u-w)\,\mathrm{d}\sigma=\widehat{U}-\widehat{W}.
$$
Then \eqref{eq:Q-representation} holds automatically.
By Lemma~\ref{lem:critical-grow-up-metric} and the second equation in \eqref{eq:system},
\[
\partial_s(D\partial_s v)=-(w-2),
\]
which yields
\[
D\partial_s v=-\int_0^s(w(\sigma,t)-2)\,\mathrm{d}\sigma=-2(\widehat{W}-s).
\]
Since
\[
\lim_{s\to0}D\partial_s v=\lim_{s\to0}D(\partial_s u-u\partial_s v)=0,
\]
we obtain
\begin{equation*}
\partial_t \widehat{U}=\frac12\int_0^s\partial_\sigma\!\left[D(\partial_\sigma u-u\partial_\sigma v)\right]\mathrm{d}\sigma
=D \partial_s^2 \widehat{U}+2(\widehat{W}-s)\partial_s \widehat{U}.
\end{equation*}
This completes the proof.
\end{proof}

The following construction is inspired by the comparison functions near the origin in \cite{2017-JEMS-TaoWinkler}.

\begin{lemma}\label{lem:critical-grow-up-subsolution}
Let
\[
c:=\frac1{20},\qquad b_0:=\frac1{1000},\qquad
a:=\frac{c}{100}=\frac1{2000},\qquad
b(t):=\frac{b_0}{1+a b_0t}.
\]
Define
\[
\underline{U}(s,t):=\frac{s}{s+b(t)h(s)},\quad
h(s):=(1-s)\bigl(1-cs(1-s)\bigr),
\]
and let $\underline{W}$ solve 
\[
\partial_t \underline{W}=\underline{U}-\underline{W},\qquad \underline{W}(s,0)=\underline{U}(s,0).
\]
Then $\partial_s \underline{U}>0$ and $\partial_t \underline{U}\ge0$ on $(0,1)\times[0,\infty)$, and
\begin{equation}\label{eq:Z-subsolution}
\partial_t \underline{U}\le D \partial_s^2 \underline{U}+2(\underline{W}-s)\partial_s \underline{U}
\quad\text{on }(0,1)\times(0,\infty).
\end{equation}
Moreover, $\underline{U}(0,t)=\underline{W}(0,t)=0$ and $\underline{U}(1,t)=\underline{W}(1,t)=1$ for $t\ge0$.
\end{lemma}

\begin{proof}
Differentiation gives
\[
b'(t)=-a b(t)^2,\quad
h^{\prime}(s)=-(1+c)+4 c s-3 c s^2,\quad
h^{\prime \prime}(s)=2 c(2-3 s).
\]
Setting $d(s,t):=s+b(t)h(s)$, we obtain
\begin{equation}\label{eq:Zt-derivative}
\partial_t \underline{U}=-\frac{sb'h}{d^2}=\frac{a b^2sh}{d^2}\geq0.
\end{equation}
\begin{equation}\label{eq:Zs-derivative}
\partial_s \underline{U}=\frac{b(h-sh')}{d^2}
=\frac{b\left[1-2cs^2(1-s)\right]}{d^2}>0.
\end{equation}
Moreover,
\begin{equation}\label{eq:Zss-derivative}
\partial_s^2 \underline{U}=-\frac{b\left\{sh''d+2(h-sh')(1+bh')\right\}}{d^3}.
\end{equation}
To prove \eqref{eq:Z-subsolution}, write
\begin{equation}\label{eq:Z-residual-decomposition}
\partial_t \underline{U}-D\partial_s^2 \underline{U}-2(\underline{W}-s)\partial_s \underline{U}=\partial_t \underline{U}+2(\underline{U}-\underline{W})\partial_s \underline{U}-R,
\end{equation}
where
\[
R:=D\partial_s^2 \underline{U}+2(\underline{U}-s)\partial_s \underline{U}.
\]
Since
\[
\underline{U}-s=\frac{s}{d}-s=\frac{s(1-s-bh)}{d},
\]
combining \eqref{eq:Zs-derivative} and \eqref{eq:Zss-derivative} gives
\begin{align}
R&=\frac{b}{d^3}\Bigl\{-Dsh''d-2D(h-sh')(1+bh')
+2s(h-sh')(1-s-bh)\Bigr\}\notag\\
&=\frac{b}{d^3}\Bigl\{-Ds^2h''
+2(h-sh')\bigl[s(1-s)-D\bigr]\notag\\
&\qquad+b\bigl[-Dsh''h-2D(h-sh')h'-2sh(h-sh')\bigr]\Bigr\}\notag\\
&=\frac{bs(1-s)}{d^3}\bigl(2sB_0+bB_1\bigr),
\label{eq:R-decomposition}
\end{align}
where
\begin{align*}
B_0&:=\tfrac12(h-sh')(1-s)
-cs\bigl(1-\tfrac12s(1-s)\bigr)(2-3s),\\
B_1&:=-s(1-s)\bigl(1-\tfrac12s(1-s)\bigr)h''
\bigl(1-cs(1-s)\bigr)\\
&\quad-2(h-sh')\left[\bigl(1-\tfrac12s(1-s)\bigr)h'
+1-cs(1-s)\right].
\end{align*}
Since $h-sh'=1-2cs^2(1-s)\ge9/10$, $7/8\le1-\tfrac12s(1-s)\le1$, and $0\le s(2-3s)\le1-s$ for $0\le s\le2/3$, we obtain
\begin{align*}
B_0&\ge\frac9{20}(1-s)-cs(2-3s)
\ge\left(\frac9{20}-c\right)(1-s)
\ge\left(\frac9{20}-c\right)\frac13=\frac2{15},
&&0\le s\le\frac23,\\
B_0&\ge\frac9{20}(1-s)
\ge\frac9{20}\cdot\frac16=\frac3{40},
&&\frac23\le s\le\frac56,\\
B_0&\ge cs\left(1-\frac12s(1-s)\right)(3s-2)
\ge c\cdot\frac56\cdot\frac78\cdot\frac12=\frac{35c}{96},
&&\frac56\le s\le1.
\end{align*}
Thus $B_0\ge c/4$.
Since $0\le s\le1$, $s(1-s)\le s$ and $4-3s\le4$, we have
\begin{align*}
&\bigl(1-\tfrac12s(1-s)\bigr)h'+1-cs(1-s)\\
&=\bigl[-1+c\bigl(s(4-3s)-1\bigr)\bigr]
\bigl(1-\tfrac12s(1-s)\bigr)+1-cs(1-s)\\
&=-c+\tfrac c2s(1-s)
+cs(4-3s)\bigl(1-\tfrac12s(1-s)\bigr)
+(\tfrac12-c)s(1-s)\\
&\le-c+\tfrac c2s+4cs+\tfrac12s
=-c+\left(\tfrac12+\tfrac{9c}{2}\right)s
=-c+\tfrac{29}{40}s,
\end{align*}
Using $h''\le4c$, $s(1-s)\le s$, $0<h-sh'\le1$, and $h-sh'\ge1-2cs$, we obtain
\begin{align*}
B_1
&\ge-4cs-2(h-sh')\left(-c+\tfrac{29}{40}s\right)\\
&=-4cs+2c(h-sh')-\tfrac{29}{20}s(h-sh')\\
&\ge-4cs+2c(1-2cs)-\tfrac{29}{20}s\\
&
=2c-\tfrac{83}{50}s
\ge2c-2s.
\end{align*}
Therefore, $b\le b_0<c/8$ implies
\[
2sB_0+bB_1
\ge\frac c2s+2cb-2bs
\ge\frac c4(s+b).
\]
Consequently, \eqref{eq:R-decomposition} yields
\begin{equation}\label{eq:R-lower-bound}
R\ge\frac c4\frac{bs(1-s)(s+b)}{d^3}.
\end{equation}
For $0\le r\le t$,
\[
d(s,r)\ge d(s,t),\qquad
\frac{b(r)}{b(t)}=1+a b(r)(t-r)\le e^{a b_0(t-r)}.
\]
Since $\partial_t \underline{W}=\underline{U}-\underline{W}$ and $\underline{W}(s,0)=\underline{U}(s,0)$,
\begin{align*}
0\le \underline{U}-\underline{W}
&=\int_0^t e^{-(t-r)}\partial_t \underline{U}(s,r)\dd r\\
&\le\frac{a b(t)^2sh}{d(s,t)^2}
\int_0^t e^{-(1-2a b_0)(t-r)}\dd r
\le\frac{2a b^2sh}{d^2}.
\end{align*}
Since $0\le h\le1-s$, $h\ge1-(1+c)s$, and $b\le b_0<1/(1+c)$, we have
\begin{equation}\label{eq:d-lower-bound}
d\ge b+\bigl[1-(1+c)b\bigr]s\ge b,
\qquad d\ge s,
\end{equation}
which yields
\[
\frac{s+b}{2}\le d\le s+b.
\]
Then, $b^2/d\le(s+b)^2/d\le2(s+b)$.
Combining \eqref{eq:Zt-derivative} and \eqref{eq:Zs-derivative} with the preceding estimate for $\underline{U}-\underline{W}$, and using $b\le b_0$, $h\le1-s$, $h-sh'\le1$ and \eqref{eq:R-lower-bound}, we obtain
\begin{align*}
\partial_t \underline{U}+2(\underline{U}-\underline{W})\partial_s \underline{U}
&\le\frac{a b^2shd}{d^3}
+\frac{4a bsh(h-sh')}{d^3}\frac{b^2}{d}\\
&\le\frac{a b_0bs(1-s)(s+b)}{d^3}
+\frac{8a bs(1-s)(s+b)}{d^3}\\
&=\left(\frac{4a b_0}{c}+\frac{32a}{c}\right)
\frac c4\frac{bs(1-s)(s+b)}{d^3}\\
&\le0.32004R<R,
\end{align*}
Substituting this into \eqref{eq:Z-residual-decomposition} yields \eqref{eq:Z-subsolution}.
Finally, $\underline{U}(0,t)=0$, $\underline{U}(1,t)=1$, and
\begin{equation}\label{eq:Y-representation}
\underline{W}(s,t)=e^{-t}\underline{U}(s,0)+\int_0^t e^{-(t-r)}\underline{U}(s,r)\dd r,
\end{equation}
yield $\underline{W}(0,t)=0$ and $\underline{W}(1,t)=1$.
\end{proof}

We now use the comparison function from Lemma~\ref{lem:critical-grow-up-subsolution} to obtain a lower bound for $\norm{u(t)}_{L^\infty}$.

\begin{proof}[Proof of Theorem~\ref{thm:critical-grow-up}]
Let $\underline{U}$ and $\underline{W}$ be defined in Lemma~\ref{lem:critical-grow-up-subsolution}, and let
\[
u_0(s):=w_0(s):=2\partial_s \underline{U}(s,0)
=\frac{2b_0[1-2cs^2(1-s)]}{[s+b_0h(s)]^2}>0.
\]
By Lemma~\ref{lem:critical-grow-up-metric}, we obtain
\[
\int_\Sigma u_0\,\mathrm{d}\mu_g=\int_\Sigma w_0\,\mathrm{d}\mu_g
=8\pi\int_0^1 \partial_s \underline{U}(s,0)\,\mathrm{d}s=8\pi.
\]
Proposition~\ref{prop:global-existence} gives a unique global classical solution, which is rotation-invariant by uniqueness.
Since $\widehat{W}(s,0)=\underline{W}(s,0)$, \eqref{eq:cumulative-system}, \eqref{eq:Q-representation}, \eqref{eq:Z-subsolution} and~\eqref{eq:Y-representation} yield
\begin{align*}
\partial_t \widehat{U}-D\partial_s^2 \widehat{U}
&=\left[2e^{-t}\widehat{W}(s,0)-2s
+2\int_0^t e^{-(t-r)}\widehat{U}(s,r)\,\mathrm{d}r\right]\partial_s \widehat{U},\\
\partial_t \underline{U}-D\partial_s^2 \underline{U}
&\le\left[2e^{-t}\widehat{W}(s,0)-2s
+2\int_0^t e^{-(t-r)}\underline{U}(s,r)\,\mathrm{d}r\right]\partial_s \underline{U}.
\end{align*}
Their initial and endpoint values satisfy
\[
\widehat{U}(s,0)=\underline{U}(s,0),\qquad
\widehat{U}(0,t)=\underline{U}(0,t)=0,\qquad \widehat{U}(1,t)=\underline{U}(1,t)=1.
\]
Fix $T<\infty$.
By \eqref{eq:d-lower-bound}, $s+bh(s)\ge b$, so
\[
0\le \partial_s \underline{U}=\frac{b(h-sh')}{[s+bh(s)]^2}\le\frac1{b(T)}
=1000+\frac{T}{2000}.
\]
Hence, the comparison principle \cite[Lemma~4.2]{2017-JEMS-TaoWinkler} yields $\widehat{U}\ge \underline{U}$ on $[0,1]\times[0,T]$.
As $T$ is arbitrary, $\widehat{U}\ge \underline{U}$ on $[0,1]\times[0,\infty)$.
Since $h(s)\le1$ for all $s\in[0,1]$ and $t\ge0$,
\[
\begin{aligned}
\|u(t)\|_{L^\infty(\Sigma)}
&\ge\frac1{b(t)}\int_0^{b(t)}u(\sigma,t)\,\mathrm{d}\sigma
=\frac{2\widehat{U}(b(t),t)}{b(t)}\\
&\ge\frac{2}{b(t)[1+h(b(t))]}
\ge\frac1{b(t)}=1000+\frac{t}{2000}.
\end{aligned}
\]
This completes the proof.
\end{proof}


\begin{thebibliography}{99}

\bibitem{1970-JTB-KellerSegel}
E.~F. Keller and L.~A. Segel,
\newblock Initiation of slime mold aggregation viewed as an instability,
\newblock \emph{J. Theoret. Biol.} \textbf{26} (1970), 399--415.

\bibitem{1992-TAMS-JaegerLuckhaus}
W.~J\"ager and S.~Luckhaus,
\newblock On explosions of solutions to a system of partial differential equations modelling chemotaxis,
\newblock \emph{Trans. Amer. Math. Soc.} \textbf{329} (1992), 819--824.

\bibitem{1995-AMSA-Nagai}
T.~Nagai,
\newblock Blow-up of radially symmetric solutions to a chemotaxis system,
\newblock \emph{Adv. Math. Sci. Appl.} \textbf{5} (1995), 581--601.

\bibitem{2001-RIMS-NagaiSenbaSuzuki}
T.~Nagai, T.~Senba, and T.~Suzuki,
\newblock Concentration behavior of blow-up solutions for a simplified
system of chemotaxis,
\newblock \emph{S\={u}rikaisekikenky\={u}sho K\={o}ky\={u}roku}
\textbf{1181} (2001), 140--176.

\bibitem{2001-JIA-Nagai}
T.~Nagai,
\newblock blow-up of nonradial solutions to parabolic-elliptic systems modeling chemotaxis in two-dimensional domains,
\newblock \emph{J. Inequal. Appl.} \textbf{6} (2001), 37--55.

\bibitem{2006-EJDE-BlanchetDolbeaultPerthame}
A.~Blanchet, J.~Dolbeault, and B.~Perthame,
\newblock Two-dimensional Keller--Segel model: optimal critical mass and qualitative properties of the solutions,
\newblock \emph{Electron. J. Differential Equations} (2006), no.~44, 1--33.

\bibitem{2018-JFA-Wei}
D.~Wei,
\newblock Global well-posedness and blow-up for the 2-D Patlak--Keller--Segel equation,
\newblock \emph{J. Funct. Anal.} \textbf{274} (2018), 388--401.

\bibitem{2001-AAM-HillenPainter}
T.~Hillen and K.~Painter,
\newblock Global existence for a parabolic chemotaxis model with prevention of overcrowding,
\newblock \emph{Adv. Appl. Math.} \textbf{26} (2001), 280--301.

\bibitem{2020-SJMA-MaheuxPierfelice}
P.~Maheux and V.~Pierfelice,
\newblock The Keller--Segel system on the two-dimensional-hyperbolic space,
\newblock \emph{SIAM J. Math. Anal.} \textbf{52} (2020), 5036--5065.

\bibitem{2023-PRSESA-Nagy}
\'{A}.~Nagy,
\newblock Stationary solutions to the Keller--Segel equation on curved planes,
\newblock \emph{Proc. Roy. Soc. Edinburgh Sect. A} \textbf{153} (2023), 327--343.

\bibitem{LiLiYan-submitted}
J.~Li, Y.~Li, and J.~Yan,
\newblock Critical mass for the Keller--Segel system on closed Riemannian surfaces,
\newblock submitted manuscript.

\bibitem{2013-BMB-StrohmTysonPowell}
S.~Strohm, R.~C. Tyson, and J.~A. Powell,
\newblock Pattern formation in a model for mountain pine beetle dispersal: linking model predictions to data,
\newblock \emph{Bull. Math. Biol.} \textbf{75} (2013), 1778--1797.

\bibitem{2017-JEMS-TaoWinkler}
Y.~Tao and M.~Winkler,
\newblock Critical mass for infinite-time aggregation in a chemotaxis model
with indirect signal production,
\newblock \emph{J. Eur. Math. Soc.} \textbf{19} (2017), 3641--3678.

\bibitem{2019-DCDSB-Laurencot}
P.~Lauren\c{c}ot,
\newblock Global bounded and unbounded solutions to a chemotaxis system with indirect signal production,
\newblock \emph{Discrete Contin. Dyn. Syst. Ser. B} \textbf{24} (2019), 6419--6444.

\bibitem{2025-JEE-XiangYang}
Z.~Xiang and L.~Yang,
\newblock Critical mass for the Cauchy problem of a chemotaxis model with indirect signal production mechanism,
\newblock \emph{J. Evol. Equ.} \textbf{25} (2025), Article~26.

\bibitem{2024-DCDSS-MaoLi}
X.~Mao and Y.~Li,
\newblock Dirac-type aggregation with full mass in a chemotaxis model,
\newblock \emph{Discrete Contin. Dyn. Syst. Ser. S} \textbf{17} (2024), 1513--1528.

\bibitem{2025-JEE-Soga}
Y.~Soga,
\newblock Concentration phenomena to a chemotaxis system with indirect signal production,
\newblock \emph{J. Evol. Equ.} \textbf{25} (2025), Article~95.

\bibitem{2003-EnseignMath-Hwang}
A.~D.~Hwang,
\newblock A symplectic look at surfaces of revolution,
\newblock \emph{Enseign. Math.} (2) \textbf{49} (2003), no.~1--2, 157--172.

\bibitem{2018-Lee}
J.~M. Lee,
\newblock \emph{Introduction to Riemannian Manifolds},
\newblock 2nd ed., Graduate Texts in Mathematics, vol.~176, Springer, Cham, 2018.

\bibitem{2007-TAMS-Yang}
Y.~Yang,
A sharp form of the Moser--Trudinger inequality on a compact Riemannian surface,
\emph{Trans. Amer. Math. Soc.} \textbf{359} (2007), 5761--5776.

\bibitem{1999-CMP-Li}
Y.~Y. Li,
\newblock Harnack type inequality: the method of moving planes,
\newblock \emph{Comm. Math. Phys.} \textbf{200} (1999), 421--444.

\bibitem{1978-BAMS-Osserman}
R.~Osserman,
\newblock The isoperimetric inequality,
\newblock \emph{Bull. Amer. Math. Soc.} \textbf{84} (1978), 1182--1238.

\bibitem{2011-Taylor}
M.~E. Taylor,
\newblock \emph{Partial Differential Equations III: Nonlinear Equations},
\newblock 2nd ed., Applied Mathematical Sciences, vol.~117,
Springer, New York, 2011.

\bibitem{2009-Grigoryan}
A.~Grigor'yan,
\newblock \emph{Heat Kernel and Analysis on Manifolds},
\newblock AMS/IP Studies in Advanced Mathematics, vol.~47,
American Mathematical Society, Providence, RI, 2009.

\bibitem{2015-TAMS-KerkyacharianPetrushev}
G.~Kerkyacharian and P.~Petrushev,
\newblock Heat kernel based decomposition of spaces of distributions in the framework of Dirichlet spaces,
\newblock \emph{Trans. Amer. Math. Soc.} \textbf{367} (2015), no.~1, 121--189.

\bibitem{1999-Folland}
G.~B. Folland,
\newblock \emph{Real Analysis: Modern Techniques and Their Applications},
\newblock 2nd ed., John Wiley \& Sons, New York, 1999.

\bibitem{2019-Lee}
D.~A. Lee,
\newblock \emph{Geometric Relativity},
\newblock Graduate Studies in Mathematics, vol.~201,
American Mathematical Society, Providence, RI, 2019.

\bibitem{1999-Hebey}
E.~Hebey,
\newblock \emph{Nonlinear Analysis on Manifolds: Sobolev Spaces and Inequalities},
\newblock Courant Lecture Notes in Mathematics, vol.~5, American Mathematical Society, Providence, RI, 1999.

\bibitem{2022-M-BehzadanHolst}
A.~Behzadan and M.~Holst,
\newblock Sobolev-Slobodeckij spaces on compact manifolds, revisited,
\newblock \emph{Mathematics} \textbf{10} (2022), no.~3, Article~522.

\bibitem{2026-arXiv-Huang}
H.~Huang,
\newblock The Cauchy problem for fully nonlinear parabolic systems on manifolds,
\newblock arXiv:1506.05030v8, 2026.

\bibitem{2006-Mueller}
R. M\"uller,
\textit{Differential Harnack Inequalities and the Ricci Flow},
EMS Series of Lectures in Mathematics,
European Mathematical Society, Z\"urich, 2006.

\bibitem{2004-BoydVandenberghe}
S.~Boyd and L.~Vandenberghe,
\newblock \emph{Convex Optimization},
\newblock Cambridge University Press, Cambridge, 2004.

\bibitem{1982-Aubin}
{
T.~Aubin,
\newblock \emph{Nonlinear Analysis on Manifolds. Monge--Amp\`ere Equations},
\newblock Grundlehren der mathematischen Wissenschaften, vol.~252,
Springer-Verlag, New York, 1982.
}

\bibitem{2012-IUMJ-GustafssonSebbar}
B.~Gustafsson and A.~Sebbar,
\newblock Critical points of Green's function and geometric function theory,
\newblock \emph{Indiana Univ. Math. J.} \textbf{61} (2012), no.~3, 939--1017.

\bibitem{2006-Jost}
J.~Jost,
\newblock \emph{Compact Riemann Surfaces: An Introduction to Contemporary Mathematics},
\newblock 3rd ed., Universitext, Springer, Berlin, 2006.

\bibitem{2015-ASNSP-Langford}
J.~J.~Langford,
\newblock Symmetrization of Poisson's equation with Neumann boundary conditions,
\newblock \emph{Ann. Sc. Norm. Super. Pisa Cl. Sci.} (5) \textbf{14} (2015), 1025--1063.

\bibitem{2001-Chavel}
I.~Chavel,
\newblock \emph{Isoperimetric Inequalities: Differential Geometric and Analytic Perspectives},
\newblock Cambridge Tracts in Mathematics, vol.~145,
Cambridge University Press, Cambridge, 2001.

\end{thebibliography}
\end{document}